\documentclass[a4paper,twoside]{article}
\usepackage{amssymb}
\usepackage{amsmath}
\usepackage[utf8]{inputenc}
\usepackage[T1]{fontenc}
\usepackage[final]{lpretex}
\usepackage[all]{xy}
\usepackage{amscd}
\usepackage{title}

\DocVersion[Tue Aug 7 2007]{1}
\DocVersion[Tue Nov 6 2007]{2}

\EndOfDump
\usepackage{htmlref}
\usepackage{url}
\newcommand\Kspol[1]{\overline{K}_0^{\mathrm{pol}}(\mathrm{Spc}_{#1})}
\newcommand\Kscheme[2][]{K_0^{#1}(\mathrm{Spc}_{#2})}
\newcommand\Kschemecompl[2][]{\widehat{K}_0^{#1}(\mathrm{Spc}_{#2})}
\newcommand\Kstack[2][]{K_0^{#1}(\mathrm{Stck}_{#2})}
\newcommand\Kcoh[1]{K_0(\mathrm{Coh}_{#1})}
\newcommand\Kcompl[1]{\widehat K_0(\mathrm{Coh}_{#1})}

\newcommand\Klat[1]{K_0(\mathrm{Ab}_{#1})}

\newcommand\Gab[1][]{L_0(\mathrm{Ab}_{#1})}
\newcommand\Gabf{L^f_0(\mathrm{Ab})}
\newcommand\NS{\mathrm{NS}}
\newcommand\CH{\mathrm{CH}}
\newcommand\Num{\mathrm{Num}}
\newcommand\Hilb{\mathrm{Hilb}}
\newcommand\sSymm{\Cal{S}\mathrm{ymm}}
\newcommand\Lclass{\mathbb{L}}
\newcommand\Bclass{\mathrm{B}}
\newcommand\Zar{\Cal{Z}\mathrm{ar}}
\newcommand\signature[1]{\mathrm{sign}_{#1}}

\newcommand\Fil{\symb{Fil}}
\newcommand\gr{\symb{gr}}
\newcommand\Bir{\symb{Bir}}
\newcommand\SBir{\symb{SBir}}
\newcommand\Conf{\symb{Conf}}
\newcommand\ArtL{\symb{ArtL}}
\newcommand\Gal{\symb{Gal}}

\let\wreath\smallint

\begin{document}

\title{A geometric invariant of a finite group} 

\author{Torsten Ekedahl} 

\address{Matematiska institutionen\\
Stockholms universitet\\ 
SE-106 91 Stockholm\\ 
Sweden}
\email{teke@math.su.se} 


\begin{abstract}
We study the class of the classifying stack of a finite group in a Grothendieck
group of algebraic stacks introduced previously. We show that this class is
trivial in a number of examples most notably for all symmetric groups. We also
give some examples where it is not trivial. The latter uses counterexamples of
Saltman and Swan to the problem of Noether.
\end{abstract}

\maketitle

In \cite{ekedahl08::class} a Grothendieck group of algebraic stacks was
introduced and was identified with a localisation of the Grothendieck group of
algebraic varieties. The purpose of this article is to investigate the classes
of stacks which can be considered to be at the opposite end of algebraic varieties,
the classifying stacks of finite groups.

After having given some general results on the class in the Grothendieck group
$\Kstack{\k}$ of algebraic stacks of the classifying stacks of finite groups we
introduce the main idea that shall be used to obtain information about such
classes. It is to use a faithful linear representation of the group that
can be exploited in an inductive approach involving the subgroups that are the
stabilisers of vectors of the representation. This is then applied to the
symmetric groups to show that the classes of their classifying stacks are equal
to $1$. In the proof a lambda structure on $\Kstack{\k}$ using ``stacky''
symmetric powers plays a crucial r\^ole. In characteristic zero we then continue
and show that that lambda structure also can be described as a certain natural
extension of the lambda structure on the Grothendieck group of varieties defined
using ordinary symmetric powers.

We go on to provide some non-triviality results. Our first result in that
direction is that the unramified Brauer group of the invariant field of a
faithful linear representation of the group (introduced in \cite{saltman84} and
identified cohomologically in \cite{bogomolov87::brauer}) is an invariant of the
class of the classifying stack and hence the examples of Saltman and Bogomolov
give examples of groups for which the class is different from $1$ (over any
field of characteristic $0$). Going backwards in time we then consider the case
of the base field being the rationals and try to fit in an example of
Swan. Using the existence of a smooth toric compactification of a (not
necessarily split) torus we show that the class in $\Kstack{\Q}$ of the
classifying stack of $\Z/47$ is non-trivial.

It is very clear from the arguments involved in both the triviality and
non-triviality that they are closely related with invariant theory for finite
groups. We finish by a short discussion on a possible actual connection.

\begin{section}{Preliminaries}

We recall some results from \scite{ekedahl08::class}. We defined a Grothendieck
group $\Kstack{\k}$ of algebraic stacks of finite type over a field $\k$ with
affine stabilisers. It is generated by isomorphism classes $\{X\}$ of such
stacks with relations $\{X\setminus Y\}=\{X\}-\{Y\}$ for a closed substack
$Y\subseteq X$ and $\{E\}=\{X\times\A^n\}$ for a vector bundle $E \to X$ of rank
$n$. It becomes a ring under $\{X\}\{Y\}=\{X\times Y\}$ and the inclusion of the
category of $\k$-schemes of finite type into that of algebraic $\k$-stacks then
induces a ring homomorphism $\Kscheme{\k} \to \Kstack{\k}$, where $\Kscheme{\k}$
is the usual Grothendieck group of $\k$-schemes. This induces (cf., \lcite[Thm.\
1.3]{ekedahl08::class}) an isomorphism $\Kscheme{\k}' \to \Kstack{\k}$, where
$\Kscheme{\k}'$ is obtained from $\Kscheme{\k}$ by inverting $\Lclass:=\{\A^1\}$
and $\Lclass^n-1$ for all $n>0$. This result is then used to get a map from
$\Kstack{\k}$ to the completion $\Kschemecompl{\k}$ of
$\Kscheme{\k}[\Lclass^{-1}]$ in the topology given by the dimension filtration.
In particular, as $\Lclass$ and $\Lclass^n-1=\Lclass^n(1-\Lclass^{-n})$ are
invertible in $\Kschemecompl{\k}$, we get a factorisation $\Kscheme{\k}\to
\Kstack{\k}\to \Kschemecompl{\k}$ of the completion map. The map $\Kstack{\k}\to
\Kschemecompl{\k}$ takes $\Fil^n\Kstack{\k}$, the subgroup spanned by stacks of
dimension $\le n$, into $\Fil^n\Kschemecompl{\k}$, the closure of the subgroup
spanned by the $\{X\}/\Lclass^m$, for schemes $X$ of dimension $\le n+m$ (cf.,
\lcite[Prop.\ 2.1]{ekedahl08::class}).

As is well-known, one can define an Euler characteristic with compact support on
$\Kschemecompl{\k}$ (see \lcite{ekedahl08::class} for a discussion on how this
is done over a general field). The recipient ring $\Kcoh{\k}$ for such an Euler
characteristic with compact support is a Grothendieck group either of mixed
Galois representations or of mixed Hodge structures. The important facts are
that we get a ring homomorphism \map{\chi_c}{\Kscheme{\k}}{\Kcoh{\k}}, that
$\Kcoh{\k}$ is a graded ring, graded by cohomological weight, and that $\chi_c$
extends to a ring homomorphism \map{\chi_c}{\Kschemecompl{\k}}{\Kcompl{\k}},
where $\Kcompl{\k}$ is the completion of $\Kcoh{\k}$ in the weight filtration,
i.e., infinite sums in the direction of negative weights are admitted.

We let $\Zar$ be the class of (connected) algebraic group schemes of finite type
all of whose torsors over any extension field of $\k$ are trivial.  Another
basic result is that if $H\subseteq G$ is a closed subgroup scheme of the finite
type subgroup scheme $G \in \Zar$, then (cf., \lcite[Prop.\ 1.1 ix]{ekedahl08::class})
\begin{equation}\label{Subgroup formula}
\{\Bclass H\}=\{G/H\}\{\Bclass G\}
\end{equation}
and when $G \in \Zar$ we have (cf., \lcite[Prop.\ 1.3]{ekedahl08::class})
\begin{equation}\label{Classifying inverse}
1=\{G\}\{\Bclass G\}.
\end{equation}
Furthermore, for any finite dimensional $\k$-algebra $L$ its group scheme of
units belongs to $\Zar$ (cf., \lcite[Prop.\ 1.3]{ekedahl08::class}). To make
later formulas explicit we shall also use the formula (cf., \lcite[Prop.\
1.1]{ekedahl08::class})
\begin{equation}
\{\GL_n\}=(\Lclass^n-\Lclass^{n-1})\cdots(\Lclass^n-1).
\end{equation}
If $G$ acts on an algebraic space, the passage to an invariant subspace commutes
with taking quotients under $G$ \emph{provided} the order of $G$ is invertible
in the base field $\k$. That is not true in general but the following
proposition shows that in a (very) special case it is so.
\begin{proposition}
\part[i]\label{subspace stratification} Let $X$ be an algebraic $\k$-space of
finite type and $Y$ a closed subspace with complement $U$. Then for every $0\le
m\le n$, the image of $Y^m\times U^{n-m}\subseteq X^n$ under the quotient map
$X^n \to X^n/\Sigma_n$ is isomorphic to $Y^m/\Sigma_m\times
U^{n-m}/\Sigma_{n-m}$.

\part[ii]\label{affine special} Let $R$ be an equivalence relation on
$\{1,2,\dots,n\}$ and let $\A^R$ be the subscheme of $\A^n$ defined by $i\sim_R
j\implies x_i=x_j$. If $\Sigma_R$ is the subgroup of $\Sigma_n$ consisting of
the permutations preserving $R$ (i.e., $i \sim_R j \implies \sigma i \sim_R
\sigma j$) then the image of $\A^R$ under the quotient map $\A^n \to
\A^n/\Sigma_R$ is isomorphic to the quotient $\A^R/\Sigma_R$.
\begin{proof}
For \DHrefpart{i} it is clear that we get an induced map $Y^m/\Sigma_m\times
U^{n-m}/\Sigma_{n-m} \to X^n/\Sigma_n$ and it will be enough to show that it is
an immersion. Now, it is clear that it is an injection on geometric points so it
is enough to show that it is unramified. Consider the induced map
\map{f_{n,m}}{X^n/\Sigma_m\times\Sigma_{n-m}}{X^n/\Sigma_n}. As the
$\Sigma_n$-stabiliser of any point of $Y^m\times U^{n-m} \subseteq X^m$ lies in
$\Sigma_m\times\Sigma_{n-m}$ we get that $f_{n,m}$ is étale along $Y^m\times
U^{n-m}$ and hence it is enough to show that $Y^m/\Sigma_m\times
U^{n-m}/\Sigma_{n-m} \to X^n/\Sigma_m\times\Sigma_{n-m}$ is unramified. This
morphism is the product of an identity map and the map $Y^m/\Sigma_m \to
X^n/\Sigma_m$. The latter map is however a closed immersion. This is well-known
and the reason is the following: The problem is local so we may assume that $X$
is affine. Then the result follows from the fact that if $U \to V$ is a
surjective map of $\k$-vector spaces, then the induced map $(U^{\tensor
n})^{\Sigma_n} \to (V^{\tensor n})^{\Sigma_n}$ is surjective and in fact a
splitting of $U \to V$ induces an equivariant splitting of $U^{\tensor n} \to
V^{\tensor n}$.

As for \DHrefpart{ii} let $\Sigma^R$ be the normal subgroup of $\Sigma_R$ fixing
$R$ (i.e., the elements that fix the $R$-equivalence classes). We start by
identifying $\A^n/\Sigma^R$ and the image of $\A^R$ in it. In the case of
characteristic zero everything is clear as then taking group quotients commutes
with passing to an (invariant) subscheme so we may assume that the base field
$\k$ has characteristic $p>0$. If $r$ is an $R$-equivalence class we put
$s_r:=|r|$ and write $s_r=p^{k_r}s'_r$ where $p\dividesnot s'_r$. We then have
that $\Sigma^R=\prod_r\Sigma_{s_r}$ (the product running over all
$R$-equivalence classes) acting diagonally on $\A^n=\prod_r\A^{s_r}$. Writing
$\A^{s_r}=\Spec\k[x_{r,i}]_{1\le i\le s_r}$ we get that $\A^n/\Sigma^R=\Spec
\k[\sigma_{r,i}]$, where the $\sigma_{r,i}$ are the elementary symmetric functions
of the $x_{r,i}$. Furthermore, the composite $\Spec\k[x_r]=\A^R \to \A^n \to
\A^n/\Sigma^R$ is induced by the graded ring homomorphism given by
$\sum_i\sigma_{r,i} \mapsto (1+x_r)^{s_r}$. This implies that the image of
$\A^R$ in $\A^n/\Sigma^R$ is equal to $\Spec \k[x_r^{p^{k_r}}]$. Now, $\Sigma_R$
acts on $\A^n/\Sigma^R$ through the quotient $\Sigma_R/\Sigma^R$. This latter
group is realised as the group of permutations of the $R$-equivalence classes
that preserve the cardinalities of the classes. Hence, as $\A^n/\Sigma^R\iso
\prod_r\A^{s_r}$ with $\Sigma_R/\Sigma^R$ permuting the factors we get that
$\A^n/\Sigma_R$ is the product of symmetric powers of the
$\A^{s_r}$. Furthermore, the image of $\A^R$ in $\A^n/\Sigma_R$ is, as we just
saw, a product of closed embeddings of $\A^1$'s in the $\A^{s_r}$. Hence we may
conclude by the fact that symmetric powers commute with closed embeddings which
was proved in \DHrefpart{i}.
\end{proof}
\end{proposition}
\end{section}
\begin{section}{The stacky and non-stacky lambda structures}

Recall that we have a lambda ring string structure on $\Kscheme{\k}$ for which
$\sigma_t(x)$ (the operation defined by $\sigma_t(x)\lambda_{-t}(x)=1$) is given
by $\sigma_t(\{X\})=\sum_{n\ge 0}\{\sigma^n(X)\}t^n$, where
$\sigma^n(X):=X^n/\Sigma_n$. (However, the references I know of assume that the
characteristic is zero. Using (\ref{subspace stratification}) that restriction
can be lifted.) We further put, for algebraic spaces $X$ and $Y$,
$\sigma^n_Y(X):=X^n\times_{\Sigma_n}\Conf^n(Y)$, where $\Conf^n(Y):=\set{(y_i)
\in Y}{\forall i\ne j\co y_i\ne y_j}$. The third part of the following
proposition is \cite[Lemma 4.4]{goettsche01::hilber} and the rest of it is
proved in essentially the same way. We shall however need the full result and
hence give the complete proof.
\begin{proposition}\label{Symmetric powers}
Let $X$ and $Y$ be algebraic $\k$-spaces of finite type and assume that $\k$ has
characteristic $0$ or $>n$.

\part[i] We have that
\begin{displaymath}
\sigma^n(\{X\}\{Y\}) = \sum_{\lambda \vdash n}\prod_i\{\sigma_Y^{n_i}(\sigma^{\lambda_i}(X))\},
\end{displaymath}
where $\lambda=[\lambda_1^{n_1},\lambda_2^{n_2},\dots,\lambda_k^{n_k}]$ with
$\lambda_1>\lambda_2>\dots>\lambda_k$ runs over the partitions of $n$.

\part[ii] There is a universal polynomial (depending only on $n$) in variables
$x_{m,\mu}$ where $\mu=(\mu_1,\dots,\mu_r)$ runs over sequences of positive integers which gives
$\sigma^n_Y(X)$ when evaluated at $x_{m,\mu}=\sigma^m(\sigma^\mu(\{X\})\{Y\})$
where $\sigma^\mu(\{X\}):=\sigma^{\mu^1}(\sigma^{\mu^2}(\dots(\{X\})))$. In
particular $\sigma_Y^n$ can be naturally extended to $\Kscheme{\k}$.

\part[iii] We have that $\sigma^n(\Lclass x)=\Lclass^n\sigma^n(x)$ for all $x \in \Kscheme{\k}$.
\begin{proof}
Given an equivalence relation $R$ on a finite set $S$
and an algebraic space $Y$ we define the subspace $Y^S_R$ of $Y^S$ by
\begin{displaymath}
Y^S_R := \set{(y_s) \in Y^S}{y_s=y_t \iff s\sim_Rt}.
\end{displaymath}
(In particular $\Conf^n(Y)=Y^n_\Delta$, where $\Delta\subseteq\{1,\dots,n\}^2$
is the diagonal.) Then $Y^n$ is the disjoint union of the $Y^n_R$ where $R$ runs
through all the equivalence relations on $\{1,\dots,n\}$ and hence $(X\times
Y)^n$ is the disjoint union of the $X^n\times Y^n_R$. The orbits under $\Sigma_n$ of
the $R$ correspond to partitions of $n$ and if let $Y^n_\lambda$ be the union of
the $Y^n_R$ over the orbit corresponding to
$\lambda=[\lambda_1^{n_1},\lambda_2^{n_2},\dots,\lambda_k^{n_k}]$, then
$(X^n\times Y^n_\lambda)/\Sigma_n$ is isomorphic to $(X^n\times Y^n_R)/N_R$, where
$R$ is the equivalence relation with equivalence classes
$\{1,\dots,\lambda_1\}$, $\{\lambda_1+1,\dots,2\lambda_1\}$,\dots,$\{(n_1-1)\lambda_1+1,\dots,n_1\lambda_1\}$
$\{n_1\lambda_1+1,\dots,n_1\lambda_1+\lambda_2\}$ etc and $N_R$ is the group of
permutation fixing $R$. We have that $N_R$ is the product
\begin{displaymath}
N_R=\left(\Sigma_{\lambda_1}\wreath\Sigma_{n_1}\right)\times\cdots\times\left(\Sigma_{\lambda_k}\wreath\Sigma_{n_k}\right)
\end{displaymath}
which implies that $(X^n\times Y^n_R)/N_R$ is isomorphic to
\begin{displaymath}
(\sigma^{\lambda_1}(X))^{n_1}\times_{\Sigma_{n_1}}\Conf^{n_1}(Y)\times
\cdots\times (\sigma^{\lambda_k}(X))^{n_k}\times_{\Sigma_{n_k}}\Conf^{n_k}(Y)
\end{displaymath}
which gives \DHrefpart{i}. Then \DHrefpart{ii} follows from \DHrefpart{i} by
induction.

Finally, applying \DHrefpart{i} to $X=\Spec\k$ gives
\begin{displaymath}
\sigma^n(\{Y\}) = \sum_{\lambda \vdash n}\prod_i\{\sigma_Y^{n_i}(1)\}
\end{displaymath}
and applied to $X=\A^1$, and using that $\sigma^m(\A^1)=\A^m$, it gives
\begin{displaymath}
\sigma^n(\Lclass\{Y\}) = \sum_{\lambda \vdash n}\prod_i\{\sigma_Y^{n_i}(\A^{\lambda_i})\}.
\end{displaymath}
Now, $\Sigma_m$ acts freely on $\Conf^m(Y)$ and acts linearly on $(\A^k)^m$ so
that $(\A^k)^m\times_{\Sigma_m}\Conf^m(Y)$ is a vector bundle over
$\Conf^m(Y)/\Sigma_m$ and hence
$\{\sigma_Y^m(\A^k)\}=\Lclass^{mk}\{\sigma_Y^m(1)\}$. This gives
\begin{displaymath}
 \sum_{\lambda \vdash n}\prod_i\{\sigma_Y^{n_i}(\A^{\lambda_i})\} =
\sum_{\lambda \vdash n}\prod_i\Lclass^{\lambda_in_i}\sigma^{n_i}_Y(1)
=\Lclass^n\sum_{\lambda \vdash n}\prod_i\{\sigma_Y^{n_i}(1)\} = \Lclass^n\sigma^n(\{Y\})
\end{displaymath}
which proves \DHrefpart{iii} in the case when $x=\{Y\}$. Now, \DHrefpart{iii}
can be expressed as the equality $\sigma_t(\Lclass x)=\sigma_{\Lclass t}(x)$ and
the multiplicativity of both sides shows that if it is true for $x=\{Y\}$, then
it is always true.
\end{proof}
\end{proposition}
\begin{remark}
\part The formula $\sigma^n(\Lclass x)=\Lclass^n\sigma^n(x)$ for all $n$ or the
equivalent $\lambda^n(\Lclass x)=\Lclass^n\lambda^n(x)$ for all $n$ would follow
from the fact that $\lambda_t(\Lclass)=1+\Lclass t$ if $\Kscheme{\k}$ were a
\emph{special} lambda ring (which it isn't).

\part The condition $p>n$ is probably necessary. We have at least that the
inverse image of a stratum of $\sigma^n(X)$ under the map $\sigma^n(X\times\A^1)
\to \sigma^n(X)$ is not necessarily a vector bundle. This is seen already for
$p=n=2$ and the ``diagonal'' stratum.
\end{remark}
We are going to use this proposition to extend this lambda structure to both
$\Kschemecompl{\k}$ and $\Kstack{\k}$. For the latter we need the following lemma.
\begin{lemma}\label{Lambda localisation}
Let $R$ be lambda ring and $r \in R$ an element such that
$\lambda_t(rx)=\lambda_{rt}(x)$ for all $x \in R$ and $\lambda_t(1)=1+t$. If $S$ is a set of integer
polynomials closed under multiplication and substitution  $t \mapsto t^n$ for
all $n>0$ and $S':=\set{f(r)}{f \in S}$, then there is a unique lambda ring
structure on the localisation $R[S'^{-1}]$ for which $R \to R':=R[S'^{-1}]$ is a
lambda ring homomorphism and for which $\lambda_t(rx)=\lambda_{rt}(x)$ for all
$x \in R'$.
\begin{proof}
Recall (cf., \cite[Exp.\ V, 2.3.2]{SGA6}) that one can define a ring structure
on $1+tA[[t]]$, functorial in the commutative ring $A$, with addition being
multiplication of power series and the multiplication being characterised (apart
from being functorial) by
being continuous in the $t$-adic topology and fulfilling
$\psi(t)\circ(1+at)=\psi(at)$ for all elements $\psi(t)$. Hence the condition
$\lambda_t(rx)=\lambda_{rt}(x)$ can be rephrased as saying that
$\lambda_t(rx)=\lambda_t(r)\circ \lambda_t(x)$. This in particular shows that
the restriction of $\lambda_t$ to the subring $U$ generated by $r$ is a ring
homomorphism and the requirement of $\lambda_t(rx)=\lambda_{rt}(x)$ for an
extension to $R'$ would make $\lambda_t$ a ring homomorphism on $U':=U[S'^{-1}]$
and the condition can then be rephrased as saying that $\lambda_t$ should be a
module homomorphism over $U'$. Hence, the uniqueness is clear and for existence
it would be enough to show that the image in $1+tR'[[t]]$ of $\lambda_t(f(r))
\in 1+tR[[t]]$ is invertible for all $f \in S$. Write $f(x)=\sum_in_ix^i$
and let $\psi(t)=1+\sum_{n>0}a_nt^n$ be a putative inverse so that we should
have $\lambda_t(f(r))\circ \psi(t)=1+t$. As $\lambda_t$ is a ring homomorphism on
$U$ this expands to
\begin{displaymath}
\prod_i\psi(r^it)^{n_i}=1+t
\end{displaymath}
and we want to show by induction over $n$ that we may choose the $a_n \in R'$ so
that this equality is valid. Now, the $t^n$-coefficient of the left hand side
equals $f(r^n)a_n+\text{(pol.\ in the }a_i, i<n, \text{ and }r)$ so that the
condition of equality for the $t^n$ coefficient may be written
\begin{displaymath}
f(r^n)a_n = \text{pol.\ in the }a_i, i<n, \text{ and }r
\end{displaymath}
and by assumption $f(r^n)$ has been inverted in $R'$.
\end{proof}
\end{lemma}
This lemma and the result of Totaro mentioned above now gives an extension of the
lambda ring structure on $\Kscheme{\k}$.
\begin{theorem}\label{Lambda extension}
The lambda ring structure on $\Kscheme{\k}$ has an extension to a lambda ring
structure on $\Kscheme{\k}'$ and $\Kscheme{\k}[\Lclass^{-1}]$ characterised by
$\lambda_t(\Lclass x)=\lambda_{\Lclass t}(x)$ for all $x \in \Kscheme{\k}'$
(resp.\ $\Kscheme{\k}[\Lclass^{-1}]$). On $\Kscheme{\k}[\Lclass^{-1}]$ it is
continuous in the dimension filtration and thus extends to a lambda ring
structure on $\Kschemecompl{\k}$.
\begin{proof}
The extensions to $\Kscheme{\k}'$ and $\Kscheme{\k}[\Lclass^{-1}]$ follow
directly from Proposition \ref{Symmetric powers} and Lemma \ref{Lambda
localisation}. One then shows that if $\dim x\le n$ for $x \in \Kscheme{\k}$,
then $\dim \lambda^i(x) \le in$ by the fact that this is obvious when $x=\{X\}$
with $\dim X \le n$ and then it follows for general elements. This and the
definition of the extension of the $\lambda^i$ to $\Kscheme{\k}[\Lclass^{-1}]$
then shows that $\dim \lambda^i(x) \le i\dim x$ remains true for
$\Kscheme{\k}[\Lclass^{-1}]$ which implies continuity.
\end{proof} 
\end{theorem}
Using the identification $\Kscheme{\k}'=\Kstack{\k}$ we in particular get a
lambda ring structure on $\Kstack{\k}$ though somewhat indirectly defined. There
is however a lambda ring structure which looks more natural. It uses the
symmetric stack powers instead of ordinary symmetric powers (the latter do not
of course make sense for algebraic stacks). Now, stack quotients of groups
acting on stacks can be somewhat tricky (particularly if the action is
non-strict). We thus start by taking by taking some time to describe the
symmetric stack powers in concrete terms.
\begin{itemize}
\item We recall the definition of the \Definition{wreath product}, $\sE\wreath
\Bclass\Sigma_n$ (which is a special case of wreath products of categories) of a
category $\sE$ and the category $\Bclass\Sigma_n$. Its objects are $n$-tuples
$(e_1,\dots,e_n)$ of objects of $\sE$. Morphisms from $(e_1,\dots,e_n)$ to
$(e'_1,\dots,e'_n)$ are tuples $(\sigma,f_1,\dots,f_n)$ where $\sigma \in
\Sigma_n$ and \map{f_i}{e_i}{e'_{\sigma^{-1}(i)}}. Composition is given by
\begin{displaymath}
(\tau,g_1,\dots,g_n)\circ(\sigma,f_1,\dots,f_n)=(\tau\sigma,g_{\sigma^{-1}(1)}\circ
f_1,\dots,g_{\sigma^{-1}(n)}\circ
f_n).
\end{displaymath}
If we have $\sE=\Bclass G$ for a group $G$ then we may also form the wreath
product of groups $G\wreath\Sigma_n$ and we have
$\Bclass(G\wreath\Sigma_n)=\Bclass G\wreath\Bclass\Sigma_n$. More generally if
the group $G$ acts on $X$ and $[X/G]$ is the action groupoid $\groupoid{G\times
X}{p_2}{m}{X}$, then $[X/G]\wreath\Bclass\Sigma_n=[X^n/(G\wreath\Sigma_n)]$
(with the obvious action of $G\wreath\Sigma_n$ on $X^n$).

\item $\sE\wreath\Bclass\Sigma_n$ clearly commutes with fibre products but not
in general with $2$-fibre products.\footnote{To us that means that $2$-fibre
squares should be commutative up to an equivalence. We shall however only use it
for groupoids where this is the same as commuting up to a natural
transformation.} Note however that there are cases when a fibre product is a
$2$-fibre product. Notably, we say that a functor \map{F}{\sA}{\sC} is a
\Definition{fibration} if it is surjective on isomorphisms. It is then easy to
show that if $\sA \to \sC$ is a fibration and $\sB \to \sC$ is any functor, then
the fibre product $\sA\times_{\sC}\sB$ is a $2$-fibre product. Now, it is
equally easy to see that if $\sA \to \sC$ is a fibration, then the induced
functor $\sA\wreath\Bclass\Sigma_n \to \sC\wreath\Bclass\Sigma_n$ is also a
fibration.

\item If $\sG$ is a pseudo-functor, then we can define a new pseudo-functor
$\sG\wreath\Bclass\Sigma_n$ given by $\left(\sG\wreath\Bclass\Sigma_n\right)(U)
= \sG(U)\wreath\Bclass\Sigma_n$. If $\sG$ is a stack then we can define the
$n$'th \Definition{symmetric stack power}, $\sSymm^n\sG$ by taking the stack
associated to this pseudo-functor. It is easy to see that if $\sG$ is a
pseudofunctor with $\sG'$ as its associated stack then the stack associated to
$\sG\wreath\Bclass\Sigma_n$ is naturally equivalent to $\sSymm^n\sG'$. (Note that
$\sSymm$ preserves $2$-fibre squares where one of the maps is just a local
fibration, i.e., locally surjective on morphisms.)

\item If $\sE$ is an algebraic stack and $\sG$ is a chart for it, i.e., a
groupoid scheme with smooth source and target maps whose associated stack is
equivalent to $\sE$ then $\sG\wreath\Bclass\Sigma_n$ is a chart for $\sSymm^n\sE$.
\end{itemize}
We start by a fairly simple result that however will be very important to us.
\begin{lemma}\label{Vector bundles}
Assume that $\sX$ is an algebraic stack and $\sE \to \sX$ a vector bundle of
rank $d$. Then $\sSymm^n\sE \to \sSymm^n\sX$ is a vector bundle of rank $nd$.
\begin{proof}
We may choose a chart $\sG=(\groupoid{\sG_1}{s}{t}{\sG_0})$ of $\sX$ such that
$\sE$ is trivial over $\sG_0$. This means that we have a functor $\sG \to
\Bclass\GL_d$ and if we let $\sG'$ be the fibre product
\begin{displaymath}
\begin{CD}
\sG' @>>> [\A^n/\GL_d]\\
@VVV   @VVV\\
\sG @>>> \Bclass\GL_d
\end{CD}
\end{displaymath}
then it is a chart for $\sE$ as $[\A^d/\GL_d] \to \Bclass\GL_d$ is a fibration
so that this is also a $2$-fibre square. Applying the wreath product we get a
cartesian square
\begin{displaymath}
\begin{CD}
\sG'\wreath\Bclass\Sigma_n @>>> [\A^d/\GL_d]=[(\A^d)^n/(\GL_d\wreath\Sigma_n)]\\
@VVV   @VVV\\
\sG\wreath\Bclass\Sigma_n @>>> \Bclass\GL_d\wreath\Bclass\Sigma_n = \Bclass(\GL_d\wreath\Sigma_n).
\end{CD}
\end{displaymath}
However, the action of $\GL_d\wreath\Sigma_n$ on $(\A^d)^n=\A^{dn}$ is linear so
that $[(\A^d)^n/(\GL_d\wreath\Sigma_n)] \to \Bclass(\GL_d\wreath\Sigma_n)$ is a
vector bundle (induced by the group homomorphism $\GL_d\wreath\Sigma_n \to
\GL_{dn}$) and therefore so is $\sSymm^n(\sE) \to \sSymm^n(\sX)$.
\end{proof}
\end{lemma}
We are now prepared for defining the more natural lambda structure on
$\Kstack{\k}$.
\begin{proposition}
There is a lambda structure $\{\lambda^n_s\}$ on $\Kstack{\k}$ with the property
that $\sigma^n_s(\{\sC\})=\{\sSymm^n\sX\}$ for any algebraic stack $\sX$. We
have that $\lambda^n_s(\Lclass x)=\Lclass^n\lambda^n_s(x)$.
\begin{proof}
We have to verify that the defining relations for $\Kstack{\k}$ are verified for
the map $\sX \mapsto \sigma^s_t(\sX):=\sum_n\{\sSymm^n\sX\}t^n$. That it is
constant on isomorphism classes is clear. If $\sE \to \sX$ is a vector bundle of
rank $d$, then by Lemma\ref{Vector bundles} we have that $\sSymm^n\sE \to
\sSymm^n\sX$ is a vector bundle of rank $dn$ so that
$\{\sSymm^n\sE\}=\{\sSymm^n\sX\times\A^{dn}\}$. If we apply this also to the
trivial vector bundle $\sX\times\A^d \to \sX$ this shows that
$\sigma^s_t(\sE)=\sigma^s_t(\sX\times\A^d)$.

Consider now the situation where $\sY \subseteq \sX$ is a closed substack of
$\sX$. We have to verify that
$\sigma^s_t(\sX)=\sigma^s_t(\sX\setminus\sY)\sigma^s_t(\sY)$. This is done in
the same way as for the lambda structure on $\Kscheme{\k}$ (cf., \cite[Lemma
3.1]{larsen04::ration}).

The last part follows from the previously established formula,
$\{\sSymm(\sX\times\A^1)\}=\{\sSymm(\sX)\times\A^n\}$.
\end{proof}
\end{proposition}
\end{section}
\begin{section}{Classifying stacks of finite group schemes}

We start by getting some formulas for the class $\Bclass G$ where $G$ is a
finite group scheme. 
\begin{proposition}\label{Formulas}
Let $V$ be an $n$-dimensional $\k$-vector space and $G \subseteq \GL(V)$ a
finite subgroup scheme.

\part $\{\Bclass G\}=\{\GL(V)/G\}/\{\GL(V)\}$.

\part The image of $\{\Bclass G\}$ in $\Kschemecompl{\k}$ is equal to
$\lim_{m\to\infty}\{V^m/G\}/\Lclass^{mn}$.

\part $\chi_c(\{\Bclass G\})=1$.
\begin{proof}
The first part is a special case of (\ref{Subgroup formula}). For the second
part we define $U_m$ to be the open subset of $V^m$ consisting of the sequences
$v_1,\dots,v_m$ for which one of the sequences $v_1,\dots,v_n$,
$v_{n+1},\dots,v_{2n}$, \dots, $v_{kn-n+1},\dots,v_{kn}$, where $k:=\lfloor
m/n\rfloor$, form a basis for $V$. Then $U_m$ is $\GL(V)$-invariant and $\GL(V)$
acts freely on it. Furthermore, the codimension of $V_m:=V^m\setminus U_m$ in
$V^m$ is equal to $k$ and hence tends to $\infty$ with $m$ and we have
\begin{multline*}
\{\Bclass G\}-\{V^m/G\}/\Lclass^{mn} = \{[V^m/G]\}/\Lclass^{mn}
-\{V^m/G\}/\Lclass^{mn}=\\
\{[U_m/G]\}/\Lclass^{mn}+\{[V'_m]\}/\Lclass^{mn}-
\{U_m/G\}/\Lclass^{mn}-\{[V_m/G]\}/\Lclass^{mn} =\\
\{[V'_m/G]\}/\Lclass^{mn}-\{V_m/G\}/\Lclass^{mn}
\stackrel{m \to \infty}{\longrightarrow}0,
\end{multline*}
where $V'_m$ is the complement of $U_m/G$ in $V^m/G$ and we have used that as
$G$ acts freely on $U_m$ we have $[U_m/G]\iso U_m/G$.

From the first part we get $\chi_c(\{\Bclass
G\})=\chi_c(\{\GL(V)/G\})/\chi_c(\{\GL(V)\})$ and by definition 
\begin{displaymath}
\chi_c(\{\GL(V)/G\})=\sum_i(-1)^i\{H^i_c(\GL(V)/G)\}
\end{displaymath}
and the same for $\GL(V)$. Now, as the coefficient for cohomology is a field of
characteristic zero we have that $H^i_c(\GL(V)/G)=H^i(\GL(V))^{G(\bk)}$. However,
the right multiplication action of $G$ on $\GL(V)$ extends to $\GL(V)$ and
$\GL(V)$ is a connected algebraic group so all elements of $\GL(V)(\bk)$ acts
trivially on $H^i_c(\GL(V))$ and hence we get $H^i_c(\GL(V)/G)=H^i_c(\GL(V))$
and in particular $\chi_c(\{\GL(V)/G\})=\chi_c(\{\GL(V)\})$.
\end{proof}
\end{proposition}
\begin{remark}
\part With small modifications of the proof and with $V^m/G$ interpreted as the
GIT-quotient $\Spec \k[V^m]^G$ the first two results remain true for arbitrary
$G$.

\part In general $\chi_c(x)$ contains a lot of information on $x$ so the fact
that $\chi_c(\{\Bclass G\})=1$ should make the equality $\{\Bclass G\}=1$ our
first guess.

\part In \cite{ekedahl08::class} a completion $\Kspol{\k}$ of $\Kscheme{\k}$ in
a stronger topology than the dimension filtration was defined and a homomorphism
$\Kstack{\k} \to \Kspol{\k}$ was shown to exist. However, the natural map
$\Kspol{\k} \to \Kschemecompl{\k}$ turned out to be injective so we don't lose
any information by computing the limit $\{V^m/G\}/\Lclass^{mn}$ in
$\Kschemecompl{\k}$ rather than in $\Kspol{\k}$. Though I haven't checked the
details it seems that this sequence actually converges in the stronger topology.
\end{remark}
The first formula is generally difficult to use as the variety $\GL(V)/G$ can be
quite complicated. There are however some cases when it (or rather(\ref{Subgroup
formula})) can be used.
\begin{proposition}\label{First triviality}
\part We have that $\{\Bclass\mu_n\}=1$ in $\Kstack{\k}$ for any field $\k$.

\part If $\k$ is a finite field and $G \in \Zar_{k}$, then $\{\Bclass G(\k)\}=1$.
\begin{proof}
We have a natural embedding $\mu_n \subseteq \mul$ and applying (\ref{Subgroup
formula}) using that $\mul \in \Zar$ we get that
$\{\Bclass\mu_n\}=\{\mul/\mu_n\}\{\Bclass\mul\}$ in $\Kstack{\k}$. Now, we have
that $\mul/\mu_n=\mul$ and by (\ref{Classifying inverse}) we have
$\{\mul\}\{\Bclass\mul\}=1\in\Kstack{\k}$.

As for the second part, the Lang torsor $G \to G$ given by $g \mapsto Fg\cdot
g^{-1}$ gives an isomorphism $G/G(\k)\iso G$ and we conclude by (\ref{Subgroup formula}).
\end{proof}
\end{proposition}
\begin{example}
The case of $\Z/n$ is more complicated. Writing $\Z/n$ as a product of its
primary components reduces to the case of $n$ being a power of a prime $p$ (as
$\Bclass(G\times H)=\Bclass(G)\times\Bclass(H)$). If $p$ divides the
characteristic we can embed $\Z/p^k$ into the truncated Witt vector scheme
$\W_m$ as the kernel of $F-1$. Hence, $\W_m/(\Z/p^k)\iso\W_m$ and as $\W_m \in
\Zar$ we get $\{\Z/p^k\}=1$. (This is a special case of Proposition \ref{First
triviality}.) We are hence left with the case when $p$ is invertible in
$\k$. Generally when $n$ is invertible in $\k$ we can consider the separable
$\k$-algebra $L:=\k[x]/(\Phi_n(x))$, where $\Phi_n(x)$ is the $n$'th cyclotomic
polynomial. We let $T$ be the torus of invertible elements in $L$. The residue
$\zeta$ of $x$ is a unit of order $n$ giving an injective map $\Z/n
\hookrightarrow T$. As $T \in \Zar$ by (\ref{Subgroup formula}) and
(\ref{Classifying inverse}) we get that
$\{\Bclass\Z/n\}=\{T/(\Z/n)\}\{T\}^{-1}$. Now, the character group of $T$, as
$\Gal(\bk/\k)$-module, is $\Z[\mu^*_n]$, where $\mu_n^*$ is the set of primitive
$n$'th roots of unity in an algebraic closure $\bk$ of $\k$. The inclusion $\Z/n
\hookrightarrow T$ is dual to the map $\Z[\mu^*_n] \to \mu_n$ which is the
inclusion $\mu^*_n\hookrightarrow\mu_n$ on generators. Hence, the character
group of $T/(\Z/n)$ is the kernel $I_n$ of this map, a Galois module first
considered by Swan (cf., \cite{swan69::invar+steen}).

\part Let us first consider only the case when $n$ is a prime $p$. Swan notices
that $I=I_p$ is a locally free module over the group ring
$\Z[\Aut(\mu_p)]\iso\Z[\Z/(p-1)]$ (this is true as the localisation equals
$\Z[\mu_p\setminus\{1\}]$ at primes different from $p$ and the order of the
Galois group is invertible at $p$). If $I$ is actually free, then $T/(\Z/p)$ is
isomorphic to $T$ and we get $\{\Bclass\Z/p\}=1$. This is certainly is true if
all rank one locally free $\Z[\Z/(p-1)]$-modules are free which in turn is true
if the class group of $\Z[\Z/(p-1)]$ is trivial. (The converse is not true
\emph{a priori} but Swan shows that for $p=47$ $I$ is not free.) Let us recall
the standard analysis of the class group of the integer group ring of a cyclic
group $C_n$ of order $n$: We have an injection $\Z[C_n] \hookrightarrow
\prod_{d|n}\Z[\zeta_d]$, where $\Z[\zeta_d]$ is the ring of integers in the
field of $d$'th roots of unity. This inclusion induces a surjection on class
groups so that in particular the class group $\Z[C_n]$ is trivial only if the
class groups of the $\Z[\zeta_d]$ are trivial. The class group of $\Z[\zeta_d]$
is trivial when $d$ is one of
$1,3,4,5,7,8,9,11,12,13,15,16,17,19,20,21,24,25,27,28,32,33,35,36,40,44,45,48,60,84$
or twice an odd number among them (as $\Z[\zeta_d]=\Z[\zeta_{2d}]$ if $d$ is
odd). Hence for an $n$ not among these, the class number of $\Z[C_n]$ is greater
than $1$. The kernel of the map on class groups is then analysed in terms of
units. We are going to show that $\Z[C_{p-1}]$ has class number $1$ for
$p=2,3,5,7,11$ and class number greater than $1$ for all other primes. For $p=2$
this is trivial and for $p=3$ it follows from the fact that the class number of
$\Z[C_q]$ is equal to that of $\Z[\zeta_q]$ for a prime $q$, a result due to Rim
(cf., \cite{rim59::modul}). Similarly the inclusion $\Z[C_4]\hookrightarrow
\Z\times\Z\times\Z[i]$ induces an isomorphism on class groups by
\cite[p.~416]{kervaire77} which takes care of $p=5$. As for $p=7,11$ we use
Milnor's Mayer-Vietoris sequence (cf.,
\cite[\S3,~Thm.~3.3]{milnor71::introd+k}): If $n$ is odd we have
$C_{2n}=C_2\times C_n$ and we get an embedding $\Z[C_{2n}] \hookrightarrow
\Z[C_n]\times\Z[C_n]$, where in the first factor we send the generator of $C_2$
to $1$ and in the second to $-1$. The subring $\Z[C_{2n}]$ of the product is
then characterised by then condition that a pair $(a,b)$ belongs to it when the
reductions of $a$ and $b$ modulo $2$ coincide. From the Milnor sequence we get
that the kernel of the map $\tilde K_0(\Z[C_{2n}]) \to \tilde K_0(\Z[C_n])^2$ is
isomorphic to the cokernel of the map $(\Z[C_n])^* \to (\Z/2[C_n])^*$. If we can
show that this cokernel is trivial for $n=3,5$ we conclude that $\Z[C_{2n}]$ has
class number $1$ again by Rim's theorem. For $n=3$ we get that
$\Z/2[C_3]=\Z/2\times\F_4$ and already $C_3 \subset (\Z[C_3])^*$ maps
surjectively onto $(\Z/2\times\F_4)^*=\F_4^*$. We have similarly for $n=5$ that
$\Z/2[C_5]=\Z/2\times\F_{16}$. Now, $(\Z/2\times\F_{16})^*=C_{15}=C_3\times C_5$
and $C_5\subset(\Z[C_5])^*$ fills out the $C_5$-part. Furthermore, $x^3+x^2-1$
is a unit in $\Z[C_5]=\Z[x]/(x^5-1)$ (with inverse $x^4+x-1$) and it maps to
$(\zeta+1)^2$ in $\F_{16}=\Z/2[x]/(x^4+x^3+x^2+x+1)$, where $\zeta$ is the
residue of $x$. As $(\zeta+1)^5=(\zeta^4+1)(\zeta+1)=\zeta^3+\zeta^2+1\ne 1$ we
see that the image of $x^3+x^2-1$ in $(\Z/2\times\F_{16})^*$ has order divisible
by $3$ and hence the map $(\Z[C_5])^* \to (\Z/2[C_5])^*$ is surjective. As for
showing that for all other $p$, the class group of $\Z[C_{p-1}]$ is non-trivial
we need only look at primes for which the class number of $\Z[\zeta_{p-1}]$ is
equal to $1$ which gives us a finite list. This list can be treated with a
suitable computer algebra system.\footnote{See,
\url{http://www.math.su.se/~teke/CyclicClassgroup.mg} for a Magma (cf.,
\cite{bosma97::magma}) script that performs this computation.} (We are only
going to use the case of triviality of $\{\Bclass\Z/p\}$.)

\part Consider now the case $n=4$. Then $I=I_4$ has a basis consisting of
$[i]-[-i]$ and $2[i]+2[-i]$ and is hence isomorphic as a Galois module to
$\Z\times\Z i$. Now, $\Z i$ is the character group of $T/\mul$ so that
$T/(\Z/4)$ is isomorphic to $\mul\times T/\mul$. As $T \to T/\mul$ is a
$\mul$-torsor we get that
$\{T/(\Z/4)\}=\{\mul\}\{T/\mul\}=\{\mul\}\{T\}\{\mul\}=\{T\}$ and hence $\{\Bclass\Z/4\}=1$.
\end{example}
\begin{example}
We assume (for simplicity) that $\k$ is algebraically closed.

\part Let $G$ be a finite subgroup of $\mul$. Then we have just seen that
$\{\Bclass G\}=1$ in $\Kstack{\k}$.

\part Let $G$ be a finite subgroup of the group of affine transformations (of
$\A^1$). We shall show by induction on the order of $G$ that $\{\Bclass G\}=1$
(in $\Kstack{\k}$ as always) and we assume that $|G|>1$. If $G$ fixes a point of
$\A^1$ then it is conjugate to a subgroup of $\mul$, a case that we have already
treated so we assume that all stabilisers of points are proper subgroup. The
stack quotient $[\A^1/G]$ is the $\A^1$-fibration associated to the universal
$G$-torsor over $\Bclass G$ so that by (\ref{Subgroup formula}) we get that
$\{[\A^1/G]\}=\Lclass\{\Bclass G\}$. Now, $\A^1$ is the disjoint union of a
finite number of $G$-orbits $O_1,\dots,O_n$ with non-trivial stabiliser and an
open subset $U$ where $G$ acts freely. This gives that
$\{[\A^1/G]\}=\sum_i\{[O_i/G]\}+[U/G]$. As $G$ acts freely on $U$ we have
$[U/G]=U/G$ and we also have $[O_i/G]\iso\Bclass G_i$ where $G_i$ is the
stabiliser of some point of $O_i$. By assumption each $G_i$ is a proper subgroup
of $G$ and hence $\{\Bclass G_i\}=1$. This gives $\{[\A^1/G]\}=n+\{U/G\}$. Now,
we have that $\A^1/G\iso \A^1$ and hence $U/G$ is open in $\A^1$ with complement
in bijection with the set of orbits $\{O_i\}$. Hence $\{U/G\}=\Lclass-n$ and we
get that $\{[\A^1/G]\}=\Lclass$. As $\Lclass$ is invertible in $\Kstack{\k}$ we
conclude from $\Lclass=\Lclass\{\Bclass G\}$ that $\{\Bclass G\}=1$.

\part Let $G$ be a finite subgroup of $\PGL_2$. We shall show by induction on
the order of $G$ that $\{\Bclass G\}=1\in\Kstack[\PGL_2]{\k}$. If $G$ has a fixed
point in its action on $\P^1$, then it is conjugate to a subgroup of the group
of affine transformations and hence we have already shown that $\{\Bclass\}=1$
so we may assume that all stabilisers are proper. As before we get that
$\{[\P^1/G]\}=(\Lclass+1)\{\Bclass G\}$ (but now only in $\Kstack[\PGL_2]{\k}$
as $\PGL_2 \notin \Zar$). We now proceed exactly as in the affine case (using
that $\P^1/G\iso\P^1$) and conclude that $\{[\P^1/G]\}=\{\P^1\}=\Lclass+1$. Now,
$\Lclass+1$, being a factor of $\Lclass^2-1$ is invertible in $\Kstack{\k}$ so
we conclude.

\part Let $G$ be the $g=1$ theta group of level a prime $p$. This means that $G$
is a non-trivial central extension of $\Z/pz$ by $\Z/p\times\Z/p$ (for $p=2$ we
get two possible groups, the dihedral and quaternion groups, depending on the
parity of the theta characteristic). Assume that $\k$ contains $p$'th root of
unity (and in particular $p\dividesnot\Char\k$). Then $G$ has an (irreducible)
representation $V$ of dimension $p$ where the central element $z$ acts by a
fixed primitive root of unity. We have that $\{[V/G]\}=\Lclass^p\{\Bclass
G\}$. On the other hand, every non-central element of $G$ has a line as fixed
point locus. Each such line is stabilised by a subgroup $\Z/p\times\Z/pz$ and
there are in all $p(p+1)$ such lines forming in all $p+1$ conjugacy classes
under the action of $G$. If $U$ is the complement of all those lines we have
that $G$ acts freely on $U$ and hence $[U/G]\iso U/G$. The complement of $U
\cup\{0\}$ is the union of the fixed point lines with the origin removed and it
divides up into conjugacy classes under $G$. If $V_L$ is such a conjugacy class
containing the line $L$ minus the origin, then
$[V_L/G]\iso[(L\setminus\{0\})/G_L]$ where $G_L$ is the stabiliser of $L$ and is
hence equal to $\Z/p\times\Z/pz$, where the first factor acts trivially and the
second freely on $L\setminus\{0\}$. Hence
\begin{displaymath}
[(L\setminus\{0\})/G_L]\iso\Bclass\Z/p\times(L\setminus\{0\})/(\Z/pz) \iso \Bclass\Z/p\times\mul
\end{displaymath}
and as $\{\Bclass\Z/p\}=1$ we get
\begin{displaymath}
\Lclass^p\{\Bclass G\} = \{U/G\}+(p+1)(\Lclass-1)+\{\Bclass G\}.
\end{displaymath}
If we also use that $V/G$ is the disjoint union of $U/G$, $p+1$ copies of $\mul$
and (the image of) $\{0\}$ we get
\begin{displaymath}
(\Lclass^p-1)\{\Bclass G\} = \{V/G\} - 1.
\end{displaymath}
As $\Lclass^p-1$ is invertible in $\Kstack{\k}$ this gives a formula for
$\{\Bclass G\}$ showing in particular that $\{\Bclass G\}=1$ precisely when
$\{V/G\}=\Lclass^p$ in $\Kscheme{\k}'$.
\end{example}
The argument used in these examples can clearly be generalised. We shall only
treat the case of the action of $\GL_n$ on $\A^n$. Hence let $G$ be a finite
group acting linearly on a finite dimensional $\k$-vector space $V$. For every
subgroup $H$ we let $V_H$ be the variety of points of $V$ whose stabiliser is
exactly $H$. It is clear that the $V_H$ for varying $H$ form a stratification of
$V$ in the strong sense of being a disjoint decomposition and the closure of
each stratum being the union of strata. More precisely, the closure of each
non-empty $V_H$ is equal to the fixed point subspace $V^H$ and
$V^H=\disjunion_{T\subseteq H}V_T$. We say that $H$ is a \Definition{stabiliser
subgroup} (wrt to $V$) if $V_H$ is non-empty. If $\sH$ is a $G$-conjugacy class
of subgroups of $G$, then we shall also denote by $V_{\sH}$ the union of the
$V_H$ for $H\in \sH$. We then have the following lemma.
\begin{lemma}
\part When $V_H$ is non-empty $N_G(H)/H$ acts freely on $V_H$ and thus
faithfully on $V^H$.

\part The stack quotient $[V/G]$ is the disjoint union of the locally closed substacks
$[V_{\sH}/G]$, where $\sH$ runs over the conjugacy classes of subgroups of $G$.

\part Each $[V_{\sH}/G]$ is isomorphic to $[V_H/N_G(H)]$ for any $H\in\sH$.

\part If $V_H$ is non-empty then $V_H$ equals $(V^H)_H$ considered as an
$N_G(H)$-representation.
\begin{proof}
These are all immediate.
\end{proof}
\end{lemma}
This lemma opens up for an inductive approach to the computation of $\{\Bclass
G\}$: We have that $\{[V/G]\}=\Lclass^{\dim V} \{\Bclass G\}$ and we can write
$[V/G]$ as the disjoint union of $[V_H/N_G(H)]$, where each $[V_H/N_G(H)]$ is
also the top stratum of the stratification of $[V^H/N_G(H)]$. Note that this
does not give a complete recursion as we are still left with the stratum
corresponding to $H=\{e\}$. If however $V$ is a faithful $G$-module, then the
action of $G$ on $V_{e}=V_{\{e\}}$ is free and hence $[V_e/G]=V_e/G$ giving at
least a recursive formula for $\{\Bclass G\}$ involving only $\{\Bclass N\}$ for
proper subgroups of $G$ and $V_e/G$. To prepare for the combinatorics that is
involved in such a formula we need to introduce some notation. Hence we define
the notion of \Definition{stabiliser flag} (with respect to the
$G$-representation $V$ and of length $n$) as follows: It is a sequence
$f=(H_0\subseteq H_1\subseteq\dots\subseteq H_n)$ of $G$ such that
\begin{itemize}
\item $H_0=\{e\}$ and

\item $H_{i+1}$ is a stabiliser subgroup of the action of $\cap_{j\le i}N_G(H_j)$ on $V^{H_i}$.
\end{itemize}
We shall also use the notation $N_G(f)$ for the normaliser of the flag (i.e.,
$\cap_{j\le n}N_G(H_j)$), $H_f$ for $H_n$, $d_f$ for $\dim V^{H_f}$ and $n_f$
for its length. Furthermore a flag is said to be \Definition{strict}
if all inclusions $H_i\subseteq H_{i+1}$ are strict.
\begin{theorem}\label{Recursion formula}
Let $G$ be a finite group and $V$ a faithful $G$-representation of dimension
$d$. Then we have that
\begin{displaymath}
\{\Bclass G\}\Lclass^d=\{V_e/G\}+\sum_{f}(-1)^{n_f}\{\Bclass N_G(f)\}\Lclass^{d_f},
\end{displaymath}
where the sum runs over $G$-conjugacy class representatives of strict stabiliser
flags of length $\ge 1$.
\begin{proof}
Note first that for an arbitrary finite group $K$ and a finite dimensional
$K$-representation $U$ over $\k$ there is a disjoint decomposition of $[V/G]$
whose pieces are the $[V_H/N_K(H)]$, where $H$ runs over the conjugacy classes
of stabiliser subgroups. This gives, together with the fact that $[V/K] \to
\Bclass K$ is a vector bundle of rank $d$, the formula
\begin{displaymath}
\{\Bclass K\}\Lclass^d=\sum_{H}\{[V_H/N_K(H)]\},
\end{displaymath}
where $H$  runs over conjugacy class representatives of stabiliser
subgroups. For  $K=G$ and $U=V$ we can separate out $H=\{e\}$ and use that $G$
acts freely on $V_e$ as $V$ is faithful to replace $[V_e/G]$ with $V_e/G$. This
gives 
\begin{displaymath}
\{\Bclass K\}\Lclass^d=\{V_e/K\}+\sum_{H}\{[V_H/N_K(H)]\},
\end{displaymath}
where the sum now is over non-trivial $H$. On the other hand we may apply the
formula with $K=N_G(f)$ and $U=V^{H_f}$ for a strict flag $f$ and writing it as
\begin{displaymath}
\{[V_{H_f}/N_G(f)]\} = \{\Bclass N_G(f)\}\Lclass^{d_f}-\sum_{f'}\{[V_{H_{f'}}/N_G(f')]\},
\end{displaymath}
where $f'$ runs over strict stabiliser flags that extend $f$ of length one more
than that of $f$. Using this formula inductively to replace the summands in the
sum gives the theorem.
\end{proof}
\end{theorem}
\begin{remark}
\part It is not true that the theorem always gives recursive formula for
$\{\Bclass G\}$ in $\Kscheme{\k}'$. The reason is that $\{\Bclass G\}$ will
appear on the right hand side, always\footnote{Except when $G=\{e\}$.} for the
flag $\{e\}\subset G$ but in general for any strict flag all of whose members
are normal subgroups of $G$. Moving all of those summands to the right we get a
formula for $\varphi(\Lclass)\{\Bclass G\}$ in terms of $\{V_e/G\}$ and
$\{\Bclass H\}$ for proper subgroups, where $\varphi$ is a monic integer
polynomial of degree $d$. The problem is that $\varphi$ may not be invertible in
$\Kscheme{\k}'$ (though there are many examples where it is). One may of course
go further and invert all monic integer polynomials in $\Lclass$.

\part It is not difficult to extend the result to finite étale group
schemes. (The only problem is that stabiliser subgroups may not be defined over
the base field. This forces one to work with finite Deligne-Mumford stacks
instead of just classifying stacks of finite groups.) For a general finite group
scheme things are different however. The proof of the theorem uses that there
are only a finite of stabiliser subgroups of $G$. For a connected finite group
scheme this is not necessarily true.
\end{remark}
The combinatorics of this formula easily become somewhat complicated. The
following lemma allows us to bypass some of these complications.
\begin{lemma}\label{Rational in L}
Two elements in $\Kstack[\sG]{\k}$ that are rational functions in $\Lclass$ are
equal when their images under $\chi_c$ are equal. In particular, for a finite
group $G$ we have that $\{\Bclass G\}=1$ in $\Kstack[\sG]{\k}$ if it is a
rational function (with integer coefficients) in $\Lclass$.
\begin{proof}
The first part follows from the fact that the ring $\Z((q^{-1}))$ of Laurent
power series in $q^{-1}=\chi_c(\Lclass^{-1})$ embeds in $\Kcompl{\k}$ and that
$\chi_c$ maps rational functions in $\Lclass$ injectively into this image. The
second part follows from the first as the rational cohomology of $G$ is trivial
we have that $\chi_c(\{\Bclass G\})=1$ in $\Kcompl{\k}$.
\end{proof}
\end{lemma}
To apply this result to the computation of $\{\Bclass G\}$ we have to fight with
the problem mentioned above. Thus for a $G$-representation $V$ we define its
\Definition{characteristic polynomial} \Definition{$\varphi_V(x)$} as 
\begin{displaymath}
x^d-\sum_f(-1)^{n_f}x^{d_f},
\end{displaymath}
where $d=\dim V$ and $f$ runs over all stabiliser flags normalised by $G$.
\begin{proposition}\label{L-polynomial}
Let $G$ be a finite group and $V$ a faithful $\k$-representation and let
$\varphi \in \Z[x]$ be a polynomial divisible by $\varphi_V$. Assume that
$\phi(\Lclass)\{\Bclass N_G(f)\}=\phi(\Lclass)$ in $\Kstack{\k}$ for all \emph{proper}
normalisers of stabiliser flags $f$. Then $\varphi(\Lclass)\{\Bclass
G\}=\varphi(\Lclass)$ precisely when
$\phi(\Lclass)\{V_e/G\}=\phi(\Lclass)\Lclass^{\dim V_e}$, where
$\phi=\varphi/\varphi_V$.
\begin{proof}
This follows immediately from Theorem \ref{Recursion formula} and Lemma
\ref{Rational in L}.
\end{proof}
\end{proposition}
Note that $N_G(f)$ in general does not act freely on $V_{H_f}$, indeed $H_f$
by definition acts trivially on it. However, $N_G(f)/H_f$ \emph{does} act freely
on $V_{H_f}$. As $H_f$ also acts trivially on $V^{H_f}$ we can apply our
recursion also to this action. This gives the following result.
\begin{proposition}\label{Quotient criterion}
Suppose that the order of $G$ is invertible in $\k$ and let $\phi \in \Z[x]$.
Assume that $\phi(\Lclass)\varphi_{V^{H_f}}(\Lclass)\{\Bclass
N_G(f)\}=\phi(\Lclass)\varphi_{V^{H_f}}(\Lclass)$ for all \emph{proper}
normalisers of stabiliser flags $f$ and $\phi(\Lclass)\{\Bclass
N_G(f)/H_f\}=\varphi(\Lclass)$ for all strict stabiliser flags $f$ of length
$\ge 1$. Then $\phi(\Lclass)\varphi_{V}(\Lclass)\{\Bclass
G\}=\phi(\Lclass)\varphi_{V}(\Lclass)$ precisely when $\phi(\Lclass)\{V/G\}$ is
a polynomial in $\Lclass$. Furthermore, this is true precisely when
$\phi(\Lclass)\{V/G\}=\phi(\Lclass)\Lclass^d$, where $d=\dim V$.
\begin{proof}
As the order of $G$ is invertible we have, for $U$ a $G$-invariant subscheme of
$V$, that the quotient map $U/G \to V/G$ maps $U/G$ isomorphically onto its
image. This means that we get a stratification of $V/G$ by the $V_H/N_G(H)$,
where $H$ runs over conjugacy class representatives of stabiliser subgroups of
$G$. We of course also have that $V_H/N_G(H)=V_H/W_H$, where we have put
$W_H:=N_G(H)/H$ so that we get $\{V/G\}=\sum_{H}\{V_H/W_H\}$. By assumption and
Proposition \ref{L-polynomial} we have that $\phi(\Lclass)\{V_H/W_H\}$ is a
polynomial in $\Lclass$ for $H\ne \{e\}$ and hence we get that
$\phi(\Lclass)\{V/G\}$ is a polynomial in $\Lclass$ precisely when
$\phi(\Lclass)\{V_e/G\}$ is. Again by Proposition \ref{L-polynomial} we conclude
that if $\phi(\Lclass)\{V/G\}$ is a polynomial in $\Lclass$ then
$\varphi(\Lclass)\{\Bclass G\}=\varphi(\Lclass)$. Finally, as the cohomology of
$V/G$ are the $G$-invariants of the cohomology of $V$ we get $V$ and $V/G$ have
the same cohomology which implies that $\chi_c(\{V/G\})=\chi_c(\Lclass^d)$ and
the last part follows from Lemma \ref{Rational in L}.
\end{proof}
\end{proposition}
\begin{subsection}{Unipotent group schemes}

We shall now show that if $G$ is a unipotent group scheme over a field $\k$,
then $\{\Bclass G\}=\Lclass^{-\dim G}$. In particular, if $G$ is finite then
$\{\Bclass G\}=1$, which fits in with the theme of the current section but the
proof of the general case is not more difficult. Recall that a group scheme is
\Definition{unipotent} if it can be embedded as a closed subgroup scheme of the
group of strictly upper triangular matrices $U_n \subseteq \GL_n$ or
equivalently every non-zero linear representation of it has a non-zero trivial
subrepresentation. We say that a smooth unipotent group scheme is
\Definition{split} if it can be obtained as a successive central extension of
$\add$'s. In particular $U_n$ is split (and over a perfect field all unipotent
group schemes are split).
\begin{proposition}\label{Affine homogeneous space}
Let $G$ be a split unipotent smooth group scheme over $\k$ and $H$ a closed
subgroup scheme. Then $G/H$ is isomorphic to $\A^{\dim G/H}$.
\begin{proof}
We prove this by induction over $\dim G$, the case of dimension $0$ being
trivial. By assumption $G$ contains a central subgroup $Z \subseteq G$
isomorphic to $\add$ and for which $G/Z$ is split. Putting $Z':=Z\cap H$ we have
a natural morphism $G/H \to (G/Z)/(H/Z')$ of homogeneous spaces. We have an
action of $Z$ on $G/H$ given by left multiplication and as $Z$ is central we get
an induced action of $Z/Z'$ making $G/H \to (G/Z)/(H/Z')$ a $Z/Z'$-torsor. Now,
a quotient of $\add$ is either $0$ or isomorphic to $\add$ again. Furthermore,
the induction assumptions imply that $(G/Z)/(H/Z')$ is isomorphic to some $\A^n$
and all $\add$-torsors over $\A^n$ are trivial so that $G/H\iso\A^{n+1}$.
\end{proof}
\end{proposition}
This immediately gives the desired result.
\begin{corollary}
Let $G$ be a unipotent group scheme over the field $\k$. Then $\{\Bclass
G\}=\Lclass^{-\dim G}$.
\begin{proof}
By assumption we can embed $G$ in some split smooth unipotent group $H$. As
$\Zar$ is closed under extensions we get that $H \in \Zar$ and therefore
$\{\Bclass G\}=\{H/G\}\{\Bclass H\}$ and $1=\{H\}\{\Bclass H\}$. Applying the
proposition to $e \subseteq H$ we get $H\iso\A^{\dim H}$ and by applying it to
$G \subseteq H$ we get $H/G\iso\A^{\dim H/G}$ and combining we get $\{\Bclass
G\}=\Lclass^{\dim(H/G)-\dim G}=\Lclass^{-\dim G}$.
\end{proof}
\end{corollary}
\begin{remark}
\part A finite group scheme $G$ is unipotent precisely when the augmentation
ideal of the dual Hopf algebra $\k[G]^*$ is nilpotent. This applies in
particular to when $G$ is an étale group scheme whose order is a power of the
characteristic of $\k$ and even more specifically when $G$ is a constant
$p$-group.

\part If $\k$ is non-perfect there are smooth connected unipotent
$1$-dimensional group schemes $G$ not isomorphic to $\add$ such as the kernel of
the morphism $\add\times\add\to \add$ given by $(x,y) \mapsto x^p-x+ty^p$, where
$t\in\k$ is not a $p$'th power. We then have that $G$ is not isomorphic as
$\k$-variety to $\A^1$ yet $\{\Bclass G\}=\Lclass^{-1}$. Of course it could
still happen that $\{G\}=\Lclass$ and hence $\{G\}\{\Bclass G\}=1$.
\end{remark}
\end{subsection}
\end{section}
\begin{section}{The symmetric groups}

We now want to apply the theorem to show that $\{\Bclass \Sigma_n\}=1$ for all
$n$. For that we shall use the standard permutation representation $V_n$ of
$\Sigma_n$. A quick way of using it to get a weaker statement is to use
Proposition \ref{Formulas} by considering $\{(V_n)^m/\Sigma_n\}$. We note that
$(V_n)^m$ is isomorphic as $\Sigma_n$-representation to $(\k^m)^n$ where the
elements of $\Sigma_n$ permutes the coordinates. Hence
$(V_n)^m/\Sigma_n\iso\sigma^n(\A^m)$ and hence by Proposition \ref{Symmetric
powers} we have $\{(V_n)^m/\Sigma_n\}=\sigma^n(\Lclass^m)=\Lclass^{mn}$ and
hence the image of $\{\Bclass\Sigma_n\}$ in $\Kcompl{\k}$ is equal to
$\lim_{n\to\infty}\Lclass^{mn}/\Lclass^{mn}=1$. Note that this argument works
only when $n!\ne 0$ in $\k$ and as we don't know if $\Kstack{\k}$ is separated
in the dimension topology we don't know how much information the equality gives
about $\{\Bclass \Sigma_n\}$ in $\Kstack{\k}$.

The idea for the full result is to use $V_n/\Sigma_n\iso\A^n$ combined with Theorem
\ref{Recursion formula}. The characteristic polynomial is $x^n-x$ which
gives an invertible element of $\Kstack{\k}$. We would then like to use
Proposition \ref{Quotient criterion} assuming inductively triviality for the
class of the classifying stack of stabiliser subgroups. This doesn't quite work
for small characteristics but there (\ref{affine special}) comes to our rescue.

The inductive nature of the proof forces us to deal with groups more general
than the symmetric groups. To prepare for this we want to interpret stabiliser
flags in more combinatorial terms. Thus if $S$ is a (finite) set, an
\Definition{equivalence flag} on $S$ is a sequence $R=(R_0\subseteq R_1\subseteq
R_2 \subseteq \dots \subseteq R_n)$ of equivalence relations on $S$ with $R_0$
being equal to the diagonal in $S\times S$. It is \Definition{strict} if the
inclusions are strict. The \Definition{flag stabiliser}, \Definition{$N_R$}, of
the flag consists of all permutations of $S$ which preserve the equivalence
relations $R_i$. We let the \Definition{residue group} of $R$ be the image of
$N_R$ in $\Sigma_{S/R_n}$ so that it is isomorphic to $N_R/H_R$, where $H_R$ is
the subgroup of $N_R$ consisting of the permutations that map each equivalence
class of $R_n$ into itself.

To the equivalence flag $R$ we associate a sequence $H_0\subseteq\dots\subseteq
H_n$ of subgroups of $\Sigma_S$, the group of permutations of $S$, namely $H_i$
is the set of permutations mapping each equivalence class of $R_i$ into
itself. We then have the following lemma.
\begin{lemma}\label{Equivalence interpretation}
The above association of a sequence of subgroups of $\Sigma_n$ to an equivalence
flag gives an $\Sigma_n$-equivariant bijection between the set of equivalence
flags on $\{1,2,\dots,n\}$ and the set of stabiliser flags of $\Sigma_n$ with
respect to its standard permutation representation. In particular the normaliser
of a stabiliser flag is the group of permutations fixing all the $R_i$.
\begin{proof}
If we think of an element of the standard permutation representation as a
function $\{1,2,\dots,n\} \to \k$ then a permutation is in the stabiliser of
such an element precisely when it takes the equivalence classes of the
equivalence relation on $\{1,2,\dots,n\}$ associated to such a function into
themselves. From this the lemma follows immediately.
\end{proof}
\end{lemma}
We can go on and give a very explicit description of flag stabilisers. For that
a \Definition{flagged set} is a set $S$ together with a sequence $R=(R_1
\subseteq R_2 \subseteq\dots\subseteq R_n)$ of equivalence relations on $S$. An
\Definition{isomorphism} of flagged sets is a bijection of the underlying sets
taking the equivalence relations of one set to the equivalence of the other (so
that in particular the lengths of the sequences of equivalence relations have to
be equal). We denote by $N_R(S)$ the automorphism group of the flagged set. We
also use $H\wreath G$ to denote the wreath product, the semi-direct product
$H^n\ltimes G$, where $G$ is a subgroup of $\Sigma_n$ and operates on $H^n$ by
permuting factors accordingly. In our cases $H$ will have a natural embedding
in $\Sigma_m$ and then $H\wreath G$ embeds naturally in $\Sigma_{mn}$.
\begin{proposition}\label{Flag stabilisers as wreath products}
Let $(S,R=(R_1 \subseteq R_2 \subseteq\dots\subseteq R_n))$ be a flagged finite
set. On each equivalence class of $R_n$ the equivalence relations $R_1 \subseteq
R_2 \subseteq\dots\subseteq R_{n-1}$ induces an equivalence flag. We define an
equivalence relation on $S$ by saying that $s \sim t$ if the induced flagged
sets $\overline{s}$ and $\overline{t}$ are isomorphic, where $\overline{s}$ and
$\overline{t}$ are the equivalence classes in $S/R_n$ containing $s$ resp.\ $t$.

\part[i] $N_R(S)$ is the product over the equivalence classes of this equivalence
relation of the flag stabilisers of the flags induced by $R$ on each such
equivalence class.

\part[ii] Assume that $S$ consists of a single equivalence class. Then $N_R(S)$ is
the wreath product $N_{R'}(S')\wreath \Sigma_t$. Here $t$ is the number of
equivalence classes of $R_n$, $S'$ is one of those equivalence classes and $R'$
is the equivalence flag of $S'$ induced by $R_1 \subseteq R_2 \subseteq\dots\subseteq R_{n-1}$.
\begin{proof}
This is essentially obvious. First an automorphism of $(S,R)$ will be able to
take an equivalence class of $R_n$ to another such equivalence class only if the
they are isomorphic as flagged sets. Hence $N_R(S)$ splits as is claimed in
\DHrefpart{i}. As for \DHrefpart{ii} we fix a flagged isomorphism between $S'$
and the other equivalence classes. This allows us to get an action of $\Sigma_t$
on $(S,R)$. On the other hand it also gives us an isomorphism between
$N_{R'}(S')^t$ and the automorphisms of $(S,R)$ taking each element of $S/R_n$
to itself.
\end{proof}
\end{proposition}
We shall say that a flagged set (possibly with an empty flag, i.e., a flag of
length $0$) $(S,R_1 \subseteq R_2 \subseteq\dots\subseteq R_n)$ is a
\Definition{homogeneous flag} if $R_n=S\times S$ and, provided that $n>1$, all
equivalence classes of $R_{n-1}$ are isomorphic as flagged sets and all the
elements of the flag are distinct and different from the least equivalence
relation (consisting only of the diagonal).  Hence, what the proposition says is
that any flag is the disjoint union of flags which become homogeneous after
removing repetitions (and possibly adding $S\times S$), that $N_R(S)$ is the
corresponding product and that for a homogeneous flag, $N_R(S)$ is an iterated
wreath product (which is just a symmetric group in the case of an empty flag).
\begin{theorem}\label{Symmetric triviality}
We have that $\{\Bclass\Sigma_n\}=1 \in \Kstack{\k}$ for all $n$.
\begin{proof}
We prove this by induction over $n$ and the case $n=1$ is trivial. We would like
to apply Proposition \ref{Quotient criterion} for the standard representation of
$\Sigma_n$. This is not possible if the characteristic is finite and $\le
n$. However, we can use (\ref{affine special}) to replace the characteristic
dependent arguments of the proposition. The characteristic polynomial for it is
$x^n-x$ and as $\Lclass^n-\Lclass$ is invertible in $\Kstack{\k}$ it is enough
to show that $\{\Bclass N_R(S)\}=1$ for all non-trivial equivalence flags $R$ on
$S$ with $|S|=n$ which we shall do by simultaneous induction over the length of
the flag and $|S|\le n$. As $\Bclass(G\times H)=\Bclass(G)\times\Bclass(H)$ we
are by Proposition \ref{Flag stabilisers as wreath products} reduced to the
case when $R$ is homogeneous. Applying again Proposition \ref{Flag stabilisers
as wreath products} we can write $N_R(S)$ as $N_{R'}(S')\wreath \Sigma_t$,
Furthermore, we have $\Bclass(N_{R'}(S')\wreath
\Sigma_t)=\Bclass(N_{R'}(S'))\wreath\Bclass\Sigma_t=\sSymm^t(\Bclass(N_{R'}(S')))$
(where we have somewhat confused constant pseudo-functors and their associated
stacks). This gives $\{\Bclass N_R(S)\}=\sigma^t_s(\{\Bclass(N_{R'}(S'))\})$
which by the induction assumption (as $R'$ is shorter than $R$) is equal to
$\sigma^t_s(1)=\{\Bclass\Sigma_t\}$ and as $R$ is non-trivial $t<n$, so by the
over-all induction assumption $\{\Bclass\Sigma_t\}=1$.
\end{proof}
\end{theorem}
This result has the following somewhat surprising corollary. For it we define,
for an algebraic space $Y$ and an algebraic stack $X$,
$\sigma_{s,Y}^n(X):=[X^n\times \Conf^n(Y)/\Sigma_n]$.
\begin{corollary}
Let $X$ be an algebraic stack and $Y$ an algebraic space of finite type both
over $\k$.

\part[i] We have that
\begin{displaymath}
\sigma_s^n(\{X\}\{Y\}) = \sum_{\lambda \vdash n}\prod_i\{\sigma_{s,Y}^{n_i}(\sigma_s^{\lambda_i}(X))\},
\end{displaymath}
where $\lambda=[\lambda_1^{n_1},\lambda_2^{n_2},\dots,\lambda_k^{n_k}]$ with
$\lambda_1>\lambda_2>\dots>\lambda_k$ runs over the partitions of $n$.

\part[ii] There is a universal polynomial (depending only on $n$) in variables
$x_{m,\mu}$ where $\mu=(\mu_1,\dots,\mu_r)$ runs over sequences of positive integers which gives
$\sigma^n_{s,Y}(X)$ when evaluated at $x_{m,\mu}=\sigma^m_s(\sigma^\mu_s(\{X\})\{Y\})$
where $\sigma_s^\mu(\{X\}):=\sigma^{\mu^1}_s(\sigma^{\mu^2}_s(\dots(\{X\})))$. In
particular $\sigma_{s,Y}^n$ can be naturally extended to $\Kscheme{\k}$.

\part[iii] The two operations $\lambda^n$ and $\lambda_s^n$ on
$\Kstack{\k}$ coincide when the  characteristic of base field is $0$ or $>n$.
\begin{proof}
For \DHrefpart{i} and \DHrefpart{ii} we follow the proof of Proposition
\ref{Symmetric powers} (as well as using its notation). Using stack quotients
instead of quotients the proof proceeds in the same fashion leading first to
\DHrefpart{i} and then using it and induction gives us \DHrefpart{ii}. Finally,
under the extra assumptions on the characteristic we can apply  \DHrefpart{i}
and \DHrefpart{ii} as well as the corresponding results of Proposition
\ref{Symmetric powers} to get
\begin{displaymath}
\sigma_s^n(\{X\}) = \sum_{\lambda \vdash
n}\prod_i\{\sigma_{s,Y}^{n_i}(\sigma_s^{\lambda_i}(1))\} =
\sum_{\lambda \vdash
n}\prod_i\{\sigma_{s,Y}^{n_i}(\sigma^{\lambda_i}(1))\} =
\sigma^n(\{X\}),
\end{displaymath}
where we have used the theorem formulated as $\sigma_s^m(1)=\sigma^m(1)=1$ for all
$m \le n$ and $\sigma_{s,Y}^m(1)=\sigma_Y^m(1)$ as the action of $\Sigma_m$ on
$\Conf^m(Y)$ is free so that $[\Conf^m(Y)/\Sigma_m]=\Conf^m(Y)/\Sigma_m$.
\end{proof}
\end{corollary}
\begin{remark}
It is probably not true that $\lambda^p=\lambda_s^p$ (on $\Kscheme{\k}$ with
values in $\Kstack{\k}$) in characteristic $p$. Indeed, the same analysis as in
the proof of the corollary shows that for a $\k$-scheme $X$ we have
$\lambda^p_s(\{X\})-\lambda^p(\{X\})=\{X\}-\{\Delta\}$, where $\Delta$ is the
image of $X$ under the composite $X \mapright{\Delta} X^p \to X^p/\Sigma_p$ and
a local computation shows that when $X$ is reduced we have that $\Delta\iso
X^{(p)}$ (where $X^{(p)}=X\times_{F_\k}\Spec\k$). Unless $X$ is defined over
$\F_p$, $X$ and $X^{(p)}$ are in general not isomorphic. For $X$ an elliptic
curve with $j$-invariant not in $\F_p$ we should probably expect $\{X\}\ne
\{X^{(p)}\}$ in $\Kstack{\k}$. (Note that if we have resolution of singularities
then an abelian variety can be recovered from its class in $\Kschemecompl{\k}$.)
\end{remark}
\end{section}
\begin{section}{Non-triviality results}

We are now going to obtain some results on the non-triviality of $\{\Bclass G\}$
for $G$ a finite group. All these results require at their basis resolution of
singularities so \emph{in this section we shall assume that $\k$ has
characteristic zero}. In \cite[Prop.\ 3.3]{ekedahl08::class} a Grothendieck
group $\Gab$ of finitely generated abelian groups with only the relations
$\{A\Dsum B\}=\{A\}+\{B\}$ was introduced. Thus $\Gab$ is the free abelian group
on the indecomposable finitely generated abelian groups $\{\Z\}$ and
$\{\Z/p^n\}$ for primes $p$ and $n>0$. For each integer $k$ a continuous
homomorphism \map{H^k(-)}{\Kschemecompl{\k}}{\Gab} was defined characterised by
the property that $H^k(\{X\}/L^{m})=\{H^{k+2m}(X,\Z)\}$, where $X$ is smooth and
proper.
\begin{remark}
This as it stands makes sense only when $\k=\C$ as only then does $H^k(X,\Z)$
exist. Otherwise, we can still define $\{H^k(X,\Z)\}$ as
$b_k(X)\{\Z\}+\sum_p\{\mathrm{tor}(H^k(X,\Z_p))\}$.
\end{remark}
We define $\Gabf$ as the subgroup of $\Gab$ spanned by the classes of finite
groups. In general any element of $\Gab$ can be written uniquely as $m\{\Z\}+a$
where $a \in \Gabf$ and we call $a$ the \Definition{torsion part} of the element
and $m\{\Z\}$ its \Definition{torsion free part} and that the element is
\Definition{torsion free} if $a=0$ and \Definition{torsion} of $m=0$.
\begin{theorem}
Let $G$ be a finite group. Then 
\begin{itemize}
\item $H^k(\{\Bclass G\})=0$ for $k>0$,

\item $H^0(\{\Bclass G\})=\{\Z\}$,

\item $H^{-1}(\{\Bclass G\})=0$,

\item $H^{-2}(\{\Bclass G\})=\{B_0(G)^\vee\}$, where $B_0(G)$ is the subgroup of
$H^2(G,\Q/\Z)$ consisting of those elements that map to zero upon restriction to
all abelian subgroups of $G$ (and $B_0(G)^\vee$ its dual as a finite group) and

\item $H^{k}(\{\Bclass G\})\in \Gabf$ for $k<0$.
\end{itemize}
\begin{proof}
Choose a faithful linear representation $V$ of $G$ of dimension $n$. By the
continuity of $H^k$ and by Proposition \ref{Formulas} we have that for large $m$
(the size depending on $k$) $H^k(\{\Bclass G\})=\{H^{k+2nm}(\{V^m/G\})\}$. Using
compactification and resolution of singularities we may write $\{V^m/G\}$ as a
linear combination $\sum_in_i\{X_i\}$ of classes of smooth and proper varieties
all of dimension $\le \dim(V^m/G)=2nm$ and then $H^k(\{\Bclass
G\})=\sum_i\{H^{k+2nm}(X_i,\Z)\}$. As $\dim X_i \le 2nm$ for all $i$ we
immediately get that $H^k(\{\Bclass G\})=0$ for $k>0$. Let $X$ be a smooth and
proper variety birational to $V^m/G$. We can then be more precise and write
$\{V^m/G\}$ as $\{X\}+\sum_in'_i\{X'_i\}$, where the $X'_i$ are smooth and
proper of dimensions $<2mn$. Hence
$\{H^{2nm}(X'_i,\Z)\}=\{H^{2nm-1}(X'_i,\Z)\}=0$ and consequently $H^0(\{\Bclass
G\})=\{H^{2nm}(X,\Z)\}=\{\Z\}$ and $H^{-1}(\{\Bclass
G\})=\{H^{2nm-1}(X,\Z)\}$. Now, $X$ being unirational is simply-connected and by
Poincaré duality $\{H^{2nm-1}(X,\Z)\}=\{H_1(X,\Z)\}=0$. Similarly,
$\{H^{2nm-2}(X'_i,\Z)\}$ is torsion free and hence the torsion part of
$H^{-2}(\{\Bclass G\})$ is equal to the torsion part of
$\{H^{2nm-2}(X,\Z)\}$. Now by Poincaré duality the torsion of
$\{H^{2nm-2}(X,\Z)\}$ is dual to the torsion of $\{H^3(X,\Z)\}$ which is the
(cohomological) Brauer group of $X$ modulo its divisible part. However, as $X$
is unirational the divisible part is zero so that the torsion of $\{H^3(X,\Z)\}$
is the Brauer group. This is a birational invariant and defined directly in
terms of the invariant field $\k(V^m)^G$ as the unramified Brauer group. This
group is isomorphic to $B_0(G)$ (cf., \cite[Thm.\ 3.1]{bogomolov87::brauer}).

Thus what remains to be shown is that $H^{k}(\{\Bclass G\})\in \Gabf$ for $k\ne
0$. Recall (cf., \cite[Prop.\ 3.2]{ekedahl08::class}) that if we compose
$H^k(-)$ with the map \map{-\Tensor\Q}{\Gab}{L_0(\Q{-}\text{vec})} given by
$\{A\}\Tensor\Q=\{A\Tensor\Q\}$ (and $L_0(\Q{-}\text{vec})$ is the Grothendieck
group of finite dimensional $\Q$-vector spaces) then for $x \in
\Kschemecompl{\k}$ we have that $(-1)^kH^k(x)\Tensor\Q$ is the weight $k$-part
of $\chi_c(x)$. (More precisely, we have a group homomorphism $\Kcoh{\k} \to
L_0(\Q{-}\text{vec})$ taking $\{V\}$ to $\dim V\{\Q\}$ and
$(-1)^kH^k(x)\Tensor\Q$ is the image of the weight $k$-part of $\chi_c(x)$ under
this map.) Now, we have (by Proposition \ref{Formulas}) that $\chi_c(\{\Bclass
G\})=1$ and hence is pure of weight $0$.
\end{proof}
\end{theorem}
\begin{corollary}
There are finite groups $G$ for which $\{\Bclass G\}\ne 1$ in $\Kschemecompl{\k}$
for all fields $\k$.
\begin{proof}
It is enough, by the theorem, to find finite groups $G$ for which $B_0(G)\ne
0$. The first such examples were given by Saltman, \cite[Thm.\ 3.6]{saltman84},
and were groups of order $p^9$ for any prime $p$ (using the definition of
$B_0(G)$ as the unramified Brauer group).  Bogomolov,
\cite{bogomolov87::brauer}, then obtained the group-cohomological description of
$B_0(G)$ and found examples of order $p^6$ with non-trivial $B_0(G)$, those
orders are minimal in terms of divisibility.
\end{proof}
\end{corollary}
\begin{subsubsection}{Abelian étale groups}

We are now going to see if we can go backwards in our previous arguments
concerning $\{\Bclass \Z/p\}$ to detect non-triviality (still in characteristic
zero of course). We can start with a somewhat more general situation: We thus
consider a finite étale $\k$-group scheme $A$ and an embedding av $A$ into a
torus $T$. Provided that $T \in \Zar_\k$ we then have $\{\Bclass
A\}=\{T/A\}/\{T\}$ so we need to obtain some information on $\{U\}$ for a torus
$U$. Using \cite{colliot-thelene05::compac+brylin+kuenn} we can find a
(projective) toric compactification $X$ of $U$ (defined over $\k$). In
particular the complement of $U$ in $X$ is a divisor with normal crossings and
we denote by $\{X_s\}_{s\in S}$ its geometric irreducible components (which are thus permuted
by the Galois group of $\k$). The fact that the Galois action in general will be
non-trivial makes for instance the inclusion-exclusion formulas a little bit
tricky and we start by introducing a formalism that will make them easier to
handle. Hence, we recall (cf., \cite{bittner04::euler}) that one defines the
Grothendieck group $\Kscheme{X}$ of schemes over a $\k$-scheme $X$ as generated
by isomorphism classes $\{Y \to X\}$ of schemes over $X$ with relations $\{Y \to
X\}=\{Y\setminus Y' \to X\}+\{Y' \to X\}$ for a closed subscheme $Y'\subseteq Y$
(we are only going to use this for $X$ a finite étale $\k$-scheme). If
\map{f}{X'}{X} is a morphism we have a pushforward map
\map{f_*}{\Kscheme{X'}}{\Kscheme{X}} given by composition with $f$ and a
pullback map \map{f^*}{\Kscheme{X}}{\Kscheme{X'}} given by pullback. The
situations we are going to encounter deal with the case when $X \to \Spec \k$ is
finite étale and $X'$ is constructed from $X$ by making some extra
choices. The push forward map will the be written as a sum over the extra
choices. This notation fulfils the expected formulas, for instance:
\begin{itemize}
\item An iterated sum can be replaced a single sum over all the choices
involved, this is transitivity $(fg)_*=f_*g_*$ of pushforward.

\item A term not depending on the extra choices may be moved out of the sum,
this is the projection formula.
\end{itemize}
As a first example, assume that $S \to \Spec\k$ is finite étale of rank $n$ and let
$\sP(S) \to \Spec\k$ be the finite étale map of flags $T_1\subset
T_2 \subset\cdots\subset T_{k} \subset S$ (of varying length) of
sub-$\k$-schemes of $S$. Then
we put
\begin{displaymath}
\signature{S} := (-1)^n\sum_{T_1\subset\dots\subset T_{k-1}\subset S}(-1)^k.
\end{displaymath}
This means more explicitly that we first define an element $(-1)^{k} \in
\Kscheme{\sP(S)}$ by writing $\sP(S)$ as the disjoint union $\sP^k(S)$ of flags
of length $k$. We then have pushforward maps induced by the
inclusions $\sP^k(S) \subseteq \sP(S)$ and we define $(-1)^{k}$ as the sum
over $k$ of the pushforwards of $(-1)^k \in \Kscheme{\sP^k(S)}$. This procedure
can be written as
\begin{displaymath}
\sum_{T_1\subset\dots\subset T_{k-1}\subset S}(-1)^{k} = \sum_{k}
\sum_{T_1\subset\dots\subset T_{k-1}\subset S}(-1)^k = \sum_{k}
(-1)^k\sum_{T_1\subset\dots\subset T_{k-1}\subset S}1.
\end{displaymath}
With this formalism we get an inclusion-exclusion formula which looks very much
like the usual one.
\begin{proposition}\label{Inclusion-exclusion}
Let $X$ be a $\k$-scheme, $S \to \Spec\k$ a finite étale map and $X'$ a closed
subscheme of $S\times X $. Let $U$ be the complement of the union $\cup_{s\in
S}X_s$, where $s$ runs over $\bk$-points of $S$ and $X_s$ is the fibre over $s$
of $X'$, then $U$ is an open $\k$-subscheme of $X$. Let $\sP(S)\to\k$ be the
étale $\k$-scheme of subschemes of $S$ with universal subscheme $T\subseteq
\sP(S)\times S$. If we put $X_T:=\cap_{s\in T}X_s$, a closed subscheme of
$\sP(S)\times X$, we have the formula
\begin{displaymath}
\{U\} = \sum_{T \subseteq S}(-1)^{|T|}\signature{T}\{X_T\} \in \Kscheme{\k}.
\end{displaymath}
\begin{proof}
Let $U_T$ be the complement in $X_T$ of the union of the $X_{T'}$ for $T\subset
T'\subseteq S$. Then $X_T$ is the disjoint union of the $U_{T'}$ for $T\subseteq
T'\subseteq S$ so that we have
\begin{displaymath}
\{X_T\} = \sum_{T\subseteq T' \subseteq S}\{U_{T'}\} =\{U_T\}+\sum_{T\subset T' \subseteq S}\{U_{T'}\}\in \Kscheme{\sP(S)}.
\end{displaymath}
The rest is then essentially to apply Möbius inversion for this formula taking
care to keep track of pushforwards. The summation formalism takes care of that
however. Thus iterating the above formula gives
\begin{multline*}
\{U_T\} = \{X_T\}-\sum_{T\subset T_1}\{U_{T_1}\} = \{X_T\}-\sum_{T\subset
T_1}\{X_{T_1}\}+\sum_{T\subset T_1\subset T_2}\{U_{T_2}\}=\dots=\\
\sum_{T\subset T_1\subset\dots\subset T_k}(-1)^k\{X_{T_k}\}=
\sum_{T\subseteq T'}\{X_{T'}\}\sum_{T\subset T_1\subset\dots\subset T_{k-1}\subset
T'}(-1)^k = \sum_{T\subseteq T'}(-1)^{|T'\setminus T|}\signature{T'\setminus T}\{X_{T'}\}
\end{multline*}
which applied to $T = \emptyset$ gives the proposition.
\end{proof}
\end{proposition}
\begin{remark}
\part In \cite{roekaeus07::some+groth} it is shown that if $S$ is of $\k$-rank
$n$, then $\signature{S}$ is equal to $(-1)^n\lambda^n(\{S\})$. We are not going
to use that but only the same formula in a Grothendieck ring of Galois
representations where we shall provide a shorter proof. In that case
$\lambda^n(\{\Z[S]\})$ is the signum character which seems conceptually
reasonable.

\part The kind of general Möbius inversion formula obtained in the proposition
will be studied more systematically elsewhere.
\end{remark}
We recall (\cite[Thm.\ 3.4]{ekedahl08::class}) that we have homomorphisms
\map{\NS^k}{\Kschemecompl{\k}}{\Gab[\k]}, where $\Gab[\k]$ is the Grothendieck
group of isomorphism classes of étale $\k$-group schemes geometrically
finitely generated and relations $\{A\Dsum B\}=\{A\}+\{B\}$. We shall denote by
$\Klat{\k}$ the Grothendieck group of the abelian category of continuous actions
of $\Gal(\bk/\k)$ on finitely generated abelian groups. We then have a
continuous group homomorphism $\Kschemecompl{\k} \to \Klat{\k}$ which is the
composite of $\NS$ and the natural map $\Gab[\k] \to \Klat{\k}$ (taking $\{A\}$
to $\{A\}$). Unsurprisingly $\Gab[\k]$ contains more information than
$\Klat{\k}$. One way of extracting more information is to consider a continuous
finite quotient $\Gal(\bk/\k) \to G$ with kernel $N$ and then regard the maps
\map{(-)^N,(-)_N}{\Gab[\k]}{\Gab[G]}, where $\Gab[G]$ is generated by
isomorphism class of finitely generated $G$-modules with the sum relation as
above, taking $\{A\}$ to the (co)invariants under $N$ $\{A^N\}$ (resp.\
$\{A_N\}$). (Note that while taking (co)invariants is not exact it is additive
which is enough to give the induced maps.) We can then compose with the natural
map $\Gab[G] \to \Klat{G}$ (with the obvious meaning of $\Klat{G}$) to get maps
\map{(-)^N,(-)_N}{\Gab[\k]}{\Klat{G}}. Note that we have a ring structure on any
$\Klat{?}$ or $\Gab[?]$ characterised by $\{A\}\{B\}=\{A\Tensor B\}$ when either
of $A$ and $B$ are torsion free and similarly a $\lambda$-ring structure
(special on $\Klat{?}$) characterised by $\lambda^i(\{A\})=\{\Lambda^i(A)\}$
when $A$ is torsion free.

In order to get an inclusion-exclusion formula we need to extend $\NS$ even
further. Hence we define, for a finite étale $\k$-scheme, $S \to \Spec\k$,
$\Gab[S]$ resp.\ $\Klat{S}$ as the Grothendieck group of étale $S$-group
schemes (or equivalently locally constant sheaves) geometrically finitely
generated (resp.\ finitely generated and torsion free) and relations $\{A\Dsum
B\}=\{A\}+\{B\}$ (resp.\ $\{C\}=\{A\}+\{B\}$ for short exact sequences
$\shex{A}{C}{B}$). The maps $\Num^k$ and $\NS^k$ then extend to this
situation. For a map \map{f}{S}{T} of finite étale $\k$-schemes the direct and
inverse images induce maps \map{f_*}{\Gab[S]}{\Gab[T]},
\map{f_*}{\Klat{S}}{\Klat{T}}, \map{f^*}{\Gab[T]}{\Gab[S]} and
\map{f^*}{\Klat{T}}{\Klat{S}}. Furthermore, the expected relations hold;
transitivity $(fg)_*=g_*f_*$ and $(fg)^*=g^*f^*$ and the projection formula
$xf_*y=f_*(f^*xy)$. Because of these relations we may (and shall) use the
summation formalism also for $\Gab[-]$ and $\Klat{-}$.
\begin{lemma}\label{Num sorites}
\part Pushforwards and pullbacks commute with $\NS^k$ and $\Num^k$.

\part If $x \in \Kscheme{\k}$ is in the subgroup generated by classes of
zero-dimensional schemes, then $\NS^k(xy)=\NS^0(x)\NS^k(y)$ for all $y \in
\Kscheme{\k}$.

\part If $S$ is a finite étale $\k$-scheme of rank $n$, then
$\Num^0(\signature{S})=\{\det\Z[S]\}$ and $\Num^k(\signature{S})=0$
for $k\ne 0$.
\begin{proof}
The first part is clear as if $\k$ is algebraically closed, pushforward and
pullback are expressed in terms of disjoint unions resp.\ direct sums, $\CH^*(-)$ and
$H^*(-,\widehat\Z)$ takes disjoint unions to direct sums and the Galois actions
are compatible.

To prove the second part we may reduce to the case when $x=\{S\}$ and in that
case $xy=f_*f^*y$, where \map{f}{S}{\Spec\k} is the structure map. Thus the
first part shows that $\NS^k(xy)=f_*f^*\NS^k(y)=f_*1\cdot\NS^k(y)$ and as
$1=\NS^0(1)$ we conclude by another application of the first part.

As for the last part, that $\Num^k(\signature{S})=0$ for $k\ne 0$ is clear and we have,
by the first part and the fact that $\Num^0(1)=1 \in \Klat{U}$ for any étale
$\k$-scheme $U$,
\begin{displaymath}
\Num^0(\signature{S})=(-1)^n\sum_{T_1\subset\dots\subset T_{k-1}\subset
S}(-1)^k\Num^0(1) = (-1)^n\sum_{T_1\subset\dots\subset T_{k-1}\subset S}(-1)^k.
\end{displaymath}
Now the sum runs over the simplices of the barycentric subdivision of the
boundary of the $n$-simplex, where we count also the empty simplex. However,
there is an extra sign so the sum is minus the equivariant Euler characteristic
of the reduced homology of the $n$-sphere. This homology is just $\Z$ in degree
$n-1$ where the symmetric group acts by the signature character.
\end{proof}
\end{lemma}
In general there is no particular relation between
these lambda ring structures, the $\NS^k$ and the lambda ring structure on
$\Kscheme{\k}$ but something can be said for special elements. Hence we define
$\ArtL_{\k} \subset \Kschemecompl{\k}$ as the closure of the ring of finite sums
$\sum_ia_i\Lclass^i$ where the $a_i$ are linear integral combinations of classes
of zero-dimensional schemes.
\begin{definition-lemma}\label{Partial lambda compatibility}
For $x \in \Kschemecompl{x}$ put $\NS_s(x) := \sum_k\NS^k(x)s^k \in
\Gab[\k]((s^{-1}))$. 

\part If $x \in \ArtL_\k$ and $y \in \Kschemecompl{\k}$ then
$\NS_s(xy)=\NS_s(x)\NS_s(y)$.

\part If we give $\Gab[\k]((s^{-1}))$ a lambda ring structure by requiring the
$\lambda^i$ to be continuous and $\lambda^i(xs^j)=\lambda^i(x)s^{ij}$ then for
$x \in \ArtL_\k$ we have $\NS_s(\lambda^i(x))=\lambda^i(\NS_s(x))$.
\begin{proof}
By continuity and additivity of the $\NS^k$ for the first part we are reduced to
the case when $x=a\Lclass^i$ and $y=b\Lclass^j$ with $a \in \ArtL_\k$ and $b \in
\Kscheme{\k}$ in which case it follows from the definition of $\NS^k$ and Lemma
\ref{Num sorites}.

For the second part assume first that the statement is true for $x$ a linear
integral combination of classes of zero-dimensional schemes. If $x =
\sum_ia_i\Lclass^i$ we have (where we extend $\NS_s$ to
$\Kschemecompl{\k}[[t]]$ by $\NS_s(\sum_ix_it^i)=\sum_i\NS_s(x_i)t^i$), using
that $\Num_s$ is a ring homomorphism on $\ArtL$ (which follows from the first part)
\begin{multline*}
\NS_s(\lambda_t(x))=\NS_s(\lambda_t(\sum_ia_i\Lclass^i))=
\NS_s(\prod_i\lambda_{\Lclass^it}(a_i)) =\\
\prod_i\lambda_{s^it}(\NS_s(a_i)) = 
\lambda_t(\sum_i\NS_s(a_i)s^i) = \lambda_t(\NS_s(x)).
\end{multline*}
What then remains to be shown is that
$\NS^0(\lambda^i(\{S\}))=\lambda^i(\NS^0(\{S\}))$ for a finite étale
$\k$-scheme $S$. As we already know that $\NS^0$ is a ring homomorphism on the
subgroup spanned by such classes it is enough to show the same formula for
$\sigma^i$ instead of $\lambda^i$. In that case
$\sigma^i(\{S\})=\{S^i/\Sigma^i\}$, $\NS^0(\{S\})=\Z[S]$,
$\NS^0(\sigma^i(\{S\}))=\Z[S^i/\Sigma_i]$ and $\sigma^i(\{\Z[S]\})=\{S^i\Z[S]\}$
and we have a natural isomorphism $S^i\Z[S]=\Z[S^i/\Sigma_i]$.
\end{proof}
\end{definition-lemma}
\begin{theorem}\label{Torus formulas}
Let $U$ be an $n$-dimensional $\k$-torus with cocharacter group $M$ considered
as a $\Gal(\bk/\k)$-module.

\part $\NS^k(\{U\})_K=(-1)^{n-k}\{\Lambda^{n-k}(M)\} \in \Klat{G}$, where $G$ is
the image of $\Gal(\bk/\k)$ in $\Aut(M)$ and $K$ is the kernel of the map
$\Gal(\bk/\k)$.

\part If $\Delta$ is an étale $\k$-sheaf of fans in $M$ which give a smooth
and proper toric compactification of $U$, $S$ the sheaf of rays of $\Delta$ and
\map{f}{\Z[S]}{M} the surjective map induced from the inclusion $S\subset M$
then $\NS^{n-1}(U)=\{\ker f\}-\{\Z[S]\} \in \Gab[\k]$.
\begin{proof}
Let $N$ be the character group of $U$ (which is the dual of $M$).
It is a fact (cf., \cite[Thm.\ p.\ 102]{fulton93::introd}) that  over $\bk$ the
cycle class map $\CH^k(Y)\Tensor\Z_p \to H^{2k}(Y,\Z_p)$ is an isomorphism for
any smooth and projective toric variety $Y$. This means that we have
$\NS^k(\{Y\})=\{\CH^k(Y)\} \in \Gab[\k]$. Here we let $\CH^k(Y)$ be the étale group
scheme with $\CH^k(Y)(\bk)=\CH^k(Y_\bk)$ with its natural Galois action. 

For the first part we shall work exclusively with $G$-modules and we shall
suppress $(-)_K$, note that we have $\CH^k(X)_K=\CH^k(X)$ and similarly for
$H^{2k}(-,\widehat\Z)$.  Now, given a smooth and proper toric $\k$-variety $Y$
corresponding to a character group $N$ with a Galois action and a fan (we think
of a fan as its set of faces) $\Delta$ in the dual of $N$ invariant under the
Galois group we have (cf., \cite[\S10]{danilov78}) that the Chow ring of $Y_\bk$
is isomorphic to $\mathrm{SR}_\Delta/(N)$, where $\mathrm{SR}_\Delta$ is the
Stanley-Reisner ring of $\Delta$ and the map $N \to \mathrm{SR}_\Delta$ is
defined as follows: By definition we have a map $\Z[S] \to \mathrm{SR}_\Delta$
into the degree $1$-part, where $S$ is the set of rays of the fan $\Delta$. The
elements of $N$ give functions of $S$ by evaluating $n \in N$ on the integral
generator of a ray. This gives a map $N \to \Z[S]^*=\Z[S]$. All these maps
respect the action of the Galois group. Furthermore, a basis for $N$ give a
regular sequence in $\mathrm{SR}_\Delta$ so that the Koszul complex for $N \to
\mathrm{SR}_\Delta$ is exact. We can then perform the standard
computation of the Hilbert series of $\mathrm{SR}_\Delta/(N)$, where one splits
up the monomial basis of it by its support, a face of $\Delta$. This gives the
following equalities of $\Klat{G}$-valued Hilbert series
\begin{equation}\label{Chow formula}
\Hilb(\CH(Y)) = \lambda_{-s}(\{N\})\Hilb(\mathrm{SR}_\Delta)=
\lambda_{-s}(\{N\})\sum_{T \in \Delta}s^{|T|}\sigma_s(\{\Z[T]\}),
\end{equation}
where $T$ runs over the faces of $\Delta$ (and as before the summation is
interpreted as a pushforward from the étale $\k$-scheme of faces of $\Delta$,
where $\Delta$ is the étale scheme of faces of $\Delta$, or equivalently
$U$-orbits on $Y$). 

We want to combine this with Proposition \ref{Inclusion-exclusion} and as the
intersection of the divisors is empty when the index set of the intersection is
not a face we do not have to sum over all subsets, only faces. Suppose thus that
$T$ is a face. Then the intersection $Y_T$ is a smooth toric variety whose torus
has character group $N_T$ which is the kernel of the surjective evaluation map
$N \to \Z[T]$ and hence using (\ref{Chow formula}) and standard formulas for
lambda rings
\begin{multline*}
\Hilb(\CH(Y_T)) = \lambda_{-s}(\{N_T\})\sum_{T\subseteq T' \in
\Delta}s^{|T'\setminus T|}\sigma_s(\{\Z[T'\setminus T]\})=\\
\lambda_{-s}(\{N\})\sigma_s(\{\Z[T]\})\sum_{T\subseteq T' \in
\Delta}s^{|T'\setminus T|}\sigma_s(\{\Z[T'\setminus T]\})=
\lambda_{-s}(N)\sum_{T\subseteq T' \in \Delta}s^{|T'\setminus
T|}\sigma_s(\{\Z[T']\}).
\end{multline*}
Using Proposition \ref{Inclusion-exclusion} we get (where we use $\signature{T}$
also for $\NS^0(\signature{T})$)
\begin{multline*}
\NS_s(U) = \lambda_{-s}(\{N\})\sum_{T\subseteq T' \in\Delta}(-1)^{|T|}\signature{T}s^{|T'\setminus
T|}\sigma_s(\{\Z[T']\}) =\\
\lambda_{-s}(\{N\}) \sum_{T\subseteq T' \in\Delta}(-1)^{|T|}\signature{T}t^{|T'\setminus
T|}\sigma_s(\{\Z[T']\}).
\end{multline*}
The sum may then be rewritten as
\begin{multline*}
\sum_{T\subseteq T' \in\Delta}(-1)^{|T|}\signature{T}s^{|T'\setminus T|}\sigma_s(\{\Z[T']\})=\\
\sum_{T' \in\Delta}(-1)^{|T'|}\signature{T'}\sigma_s(\{\Z[T']\})\sum_{T\subseteq
T'}\signature{T'\setminus T}(-s)^{|T'\setminus T|},
\end{multline*}
where we have used that $\signature{T'}=\signature{T'}\signature{T'\setminus T}$
and $\signature{T'}^2=1$. Now,
\begin{displaymath}
\sum_{T\subseteq T'}\signature{T'\setminus
T}(-s)^{|T'\setminus T|}=
\sum_{0\le k\le m}(-s)^{m-k}\sum_{\stackrel{T\subseteq
T'}{|T|=m-k}}\signature{T}=
\sum_{0\le k\le m}(-s)^{k}\sum_{\stackrel{T\subseteq
T'}{|T|=k}}\signature{T},
\end{displaymath}
where $m:=|T'|$, and the inner sum is just $\{\Lambda^k\Z[T']\}$ so that the
full sum is equal to $\lambda_{-s}(\{\Z[T']\})$. This gives
\begin{multline*}
\Num_s(U) =\lambda_{-s}(\{N\}) \sum_{T' \in
\Delta}(-1)^{|T'|}\signature{T'}\sigma_s(\{\Z[T']\})\lambda_{-s}(\{\Z[T']\})=\\
\lambda_{-s}(\{N\}) \sum_{T' \in
\Delta}(-1)^{|T'|}\signature{T'}.
\end{multline*}
The sum is the reduced Euler characteristic of $\Delta$ (again with an extra sign) considered as a
simplicial complex (it is unoriented which is the reason why one has to insert
the $\signature{T'}$). Now, $\Delta$ is a triangulation of the sphere of
directions of $M\Tensor\R$ so that its reduced homology is concentrated in
degree $n-1$ and is $\Z$ there. In general if $V$ is an $n$-dimensional real
vector space and $S$ its sphere of directions (i.e., $V\setminus\{0\}$ modulo
rays) then there is a canonical identification
$H_{n-1}(S,\R)=\Lambda^nV$. Hence, in our case the Galois action on
$H_{n-1}(\Delta,\R)$ is the same as the action on
$\Lambda^n(M)\Tensor\R$. However, as the action on both $H_{n-1}(\Delta,\Z)$ and
$\Lambda^n(M)$ are given by characters this implies that $H_{n-1}(\Delta,\Z)\iso
\Lambda^n(M)$ and thus
\begin{displaymath}
 \sum_{T' \in
\Delta}(-1)^{T'}\signature{T'} = (-1)^n\lambda^n(\{M\})
\end{displaymath}
giving
\begin{displaymath}
\Num_s(U) = \lambda_{-s}(\{N\})\lambda^n(\{M\})(-1)^n
\end{displaymath}
and separating powers of $t$ gives
\begin{displaymath}
\Num^k(\{U\}) = (-1)^{n+k-1}\{\Lambda^kN\Tensor\Lambda^nM\}.
\end{displaymath}
Now, the wedge product gives isomorphisms
$\Lambda^{n-k}M=(\Lambda^{k}M)^\vee\Tensor\Lambda^nM$ and as we also have natural
isomorphisms $(\Lambda^{k}M)^\vee=\Lambda^{k}M^\vee=\Lambda^{k}N$ we get
\begin{displaymath}
\Num^k(\{U\}) = (-1)^{n+k}\{\Lambda^{n-k}(M)\} =  (-1)^{n-k}\{\Lambda^{n-k}(M)\}
\end{displaymath}
which proves the first part.

As for the second part we choose again a toric compactification $X$ of $U$ and
apply $\NS^k$ to the formula of Proposition \ref{Inclusion-exclusion} for
$k=n-1$. Using that $\CH^{n-1}(Y)=0$ if $\dim Y<n-1$ and $\CH^{n-1}(Y)=\Z[S]$ if
$\dim Y=n-1$ and $S$ the $\k$-scheme of $n-1$-dimensional irreducible
geometric components we get that 
\begin{displaymath}
\NS^{n-1}(\{U\}) = \{\CH^{n-1}(X)\}-\{\Z[S]\},
\end{displaymath}
where $S$ is the scheme of irreducible geometric components of $X\setminus U$
(or equivalently the rays of the fan $\Delta$ of $X$). Now, multiplication in
the Chow ring gives a perfect duality (as the Chow groups tensored with
$\widehat\Z$ is the cohomology of $X$) between $\CH^{n-1}(X)$ and
$\CH^1(X)=\Pic(X)$. As a split torus has trivial Picard group we get an exact
sequence
\begin{displaymath}
\shex{N}{\Z[S]}{\Pic(X)}
\end{displaymath}
and dualising gives an exact sequence
\begin{displaymath}
\shex{\CH^{n-1}(X)}{\Z[S]}{M}
\end{displaymath}
which finishes the proof.
\end{proof}
\end{theorem}
To be able to formulate the next theorem let us introduce an involution
$(-)^\vee$ of $\Klat{?}$ given by $\{A\}^\vee =
\{\Hom_\Z(A,\Z)\}-\{\Ext^1_\Z(A,\Z)\}$. In particular, if $A$ is finite we have
that $\{A\}^\vee = -\{\Hom(A,\Q/\Z)\}$, i.e., minus the usual dual of $A$.
\begin{theorem}
Let $A$ be a finite étale commutative $\k$-group scheme. If $G$ is the image
of $\Gal(\bk/\k)$ in $\Aut(A)$ then $\NS_s(\{\Bclass
A\})=\lambda_{-s^{-1}}(\{A\})\in\Klat{G}$.
\begin{proof}
We shall consistently confuse étale commutative $\k$-group schemes with
representations of $\Gal(\bk/\k)$ (and generally étale $k$-schemes with
$\Gal(\bk/\k)$-sets). Further, as we shall only consider such sets for which the
action factors through $G$ we shall also confuse with $G$-representations
($G$-sets). Pick a $\k$-subscheme $S$ of $A'$ such that it generates $A'$ (we
may choose $S=A'$). This gives us a surjective homomorphism of étale sheaves
$N:=\Z[S] \to A'$ and let $N'$ be its kernel. Let $T$ resp.\ $T'$ be the
$\k$-tori corresponding to $N$ resp.\ $N'$. In particular $T$ is the torus of
units in $\Gamma(S,\sO_S)$ so that $T \in \Zar$. The exact sequence
$\shex{N'}{N}{A'}$ gives an exact sequence of group schemes $\shex{A}{T}{T'}$
and hence, as $T \in \Zar$, we have $\{\Bclass A\}=\{T'\}\{T\}^{-1}$. If we let
$M$ resp.\ $M'$ be the cocharacter groups of $T$ resp.\ $T'$ (i.e., the duals of
$N$ resp.\ $N'$) we get from Theorem \ref{Torus formulas} that
$\NS_s(\{T\})=s^n\lambda_{-s^{-1}}(\{M\})$ (resp.\
$\NS_s(\{T'\})=s^n\lambda_{-s^{-1}}(\{M'\})$) where $n:=\dim T=\dim T'$. In
\cite{roekaeus07::some+groth} a formula for $\{T\} \in \Kscheme{\k}$ is proven
which implies that $\{T\} \in \ArtL_\k$. Hence by Lemma \ref{Partial lambda
compatibility} we have
\begin{displaymath}
\NS_s(\{\Bclass A\})=\NS_s(\{T'\})\NS_s(\{T\})^{-1}=\lambda_{-s^{-1}}(\{M'\})\lambda_{-s^{-1}}(\{M\})^{-1}
\end{displaymath}
and as dualising the exact sequence $\shex{N'}{N}{A'}$ gives us an exact
sequence $\shex{M}{M'}{A}$ we have
\begin{displaymath}
\lambda_{-s^{-1}}(\{M'\})\lambda_{-s^{-1}}(\{M\})^{-1}= 
\lambda_{-s^{-1}}(\{M'\}-\{M\})=\lambda_{-s^{-1}}(\{A\}).
\end{displaymath}
\end{proof}
\end{theorem}
From this theorem it follows that there are many étale group schemes $A$ (over
suitable fields $\k$) for which $\{\Bclass A\}\ne 1$. In fact, by the
localisation sequence of algebraic $K$-theory, for any finite group $G$ the
kernel of the map $\Klat{\Z[G]} \to \Klat{\Q[G]}$ is generated by classes of
torsion $G$-modules and this kernel is often non-trivial. Then for any torsion
$G$-module $A$ whose class $\{A\} \in \Klat{G}$ is non-trivial and any field
$\k$ whose Galois group has $G$ as a quotient we get an étale group scheme $A$
for which $\{\Bclass A\}$ is different from $1$. We content ourselves by noting
that Swan's example does indeed give a non-triviality result also in our
context.
\begin{corollary}
$\{\Bclass\Z/47\Z\} \ne 1$ in $\Kschemecompl{\Q}$.
\begin{proof}
(The calculations to follow are of course standard but are given for lack of a
reference.) If $\{\Bclass\Z/47\Z\}=1$, then by the theorem
$\lambda_t(\{\mu_{47}\})=1$ in $\Klat{(\Z/47\Z)^*}$ and in particular
$\{\mu_{47}\}=0$. Now, the maximal order of $\Z[\Z/23]$ is
$\Z\times\Z[\zeta_{23}]$, where $\zeta_{23}$ is a primitive $23$'rd root of
unity. Furthermore, the map $K_0(\Z\times\Z[\zeta_{23}]) \to K_0(\Z[\Z/23])$
induced by the inclusion is an isomorphism (cf., \cite[Cor.\
39:26]{curtis87::method}). The action of $\Z[\Z/23]$ factors through the
inclusion $\Z[\Z/23]\subseteq \Z\times\Z[\zeta_{23}]$ (by direct inspection or
as the index of that inclusion is relatively prime to $47$) and by the
isomorphism above it is enough to show that $\{\mu_{47}\}\ne 0\in
K_0(\Z\times\Z[\zeta_{23}])=K_0(\Z[\zeta_{23}])$. Now, picking a generator
$\zeta$ for $\mu_{47}$ we get a surjective $\Z/23$-equivariant map
$\Z[\zeta_{23}]\to \mu_{47}$ mapping $\zeta_{23}^i$ to $\zeta^{ia}$, where $a
\in (\Z/47)^*$ is an element of order $23$, and its kernel is an ideal $I$ of
$\Z[\zeta_{23}]$ so that $\{\mu_{47}\} = \{\Z[\zeta_{23}]\}-\{I\}$ and as
$K_0(\Z[\zeta_{23}])=\Z\Dsum\Pic(\Z[\zeta_{23}])$ so that $\{\mu_{47}\} \ne 0$
precisely when $I$ is non-principal. That it is non-principal is \cite[Thm.\
3]{swan69::invar+steen}.
\end{proof}
\end{corollary}
\begin{remark}
\part The results on the classes in $\Klat{G}$ are somewhat disappointing from
the point of view of finding invariants as they say in particular that the
values of the $\NS^k$ are determined by the \emph{class} of the character group
of the torus resp.\ finite group scheme. They must of course be determined by
the character group itself but the higher invariants could possibly contain more
information than the first one (which is essentially the class of the character
group). 

In the case of the invariants in $\Gab[\k]$ the situation seems to be different
however. To begin with for an $n$-dimensional torus $U$ with cocharacter group
$M$ it is not necessarily the case that $\NS^{n-1}(\{U\})=-\{M\}$ in
$\Gab[\k]$. Indeed, consider a field for which we have a quotient $\Gal(\bk/\k)
\to \Sigma_n$ with kernel $K$ and let
$M:=\Z[\{1,2,\dots,n\}]/\Z([1]+\dots+[n])$. Here we can compactify the
corresponding torus by using the standard $\Sigma_n$-action on $\P^{n-1}$ and
twisting by the $\Sigma_n$-torsor over $\k$ given by the surjection above. In
other words, the generators of the rays of the fan are exactly the residues of
the $[i]$ (and the cones are generated by the images of $\{1,2,\dots,n\}$ minus
one element). Then Theorem \ref{Torus formulas} gives that $\NS^{n-2}(U)$ is
equal to $\{\Z\}-\{\Z[\{1,2,\dots,n\}]\}$. This is not equal $-\{M\}$ in
$\Gab[\k]$. Indeed, if that were the case we could apply
\map{\Hom_\Z(-,\Z)}{\Gab[\k]}{\Gab[\k]} to the equality and get that
$\{N\}=\{\Z[\{1,2,\dots,n\}]\}-\{\Z\}$. Then applying $(-)_K$ shows that we
would have the same equality in $\Gab[\Sigma_n]$. However, we have a map
\map{H^1(\Sigma_n,-)}{\Gab[\Sigma_n]}{\Gab} and the right hand side maps to zero
whereas $H^1(\Sigma_n,N)\ne 0$. Indeed, we have an exact sequence
$\shex{\Z}{\Z[\{1,2,\dots,n\}]}{N}$ and the long exact sequence of cohomology
shows that $H^1(\Sigma_n,N)=\Z/n$. In this case however we have that
$\NS^k(\{U\})=(-1)^{n-1-k}\lambda^{n-1-k}(\NS^{n-2}(\{U\}))$ in $\Gab[\k]$. This
follows from a formula (see \cite[Rmk.\ 4.4]{roekaeus07::some+groth}) for
$\{U\}$ in $\Kscheme{\k}$. I doubt that such a formula would be true in general.

\part It seems almost miraculous that the original terms simplify to the actual
results. It would be nice to have a more \emph{a priori} proof of the final
result. It is clear that one would get the same result if instead of working
over a non-algebraically closed field, one worked $G$-equivariantly over the
complex numbers so that we can use integral cohomology. We then have a
$G$-equivariant spectral sequence
\begin{displaymath}
E_1^{i,j}=\prod_{|T|=i} H^j(X_T,\Z) \implies H^*_c(U)
\end{displaymath}
coming from covering the complement of $U$ by the $X_s$ (and the product is over
all subsets of $S$ including the empty one). (Note that there is an implicit
signature character present in the $G$-action on the $E_1$ as can be seen
already in Meyer-Vietoris sequences.) This spectral sequence degenerates
rationally at the $E_2$-level for weight reasons and as each $H^{2n-k}_c(U)$ is
pure of weight $2n-2k$ (for $k\ge n$) only one row contributes to a given
$H^{2n-k}_c(U)$. Provided that the $E_2$-term is torsion free this remains true
over the integers and would immediately give the theorem. Torsion freeness seems
plausible; the zeroth row for instance is just the reduced chain complex of the
fan $\Delta$ which is a triangulation of the sphere.
\end{remark}
\end{subsubsection}
\end{section}
\begin{section}{Relations with invariant theory}
\
We have seen that the problem on whether $\{\Bclass G\}=1$ for a finite group
$G$ is closely related to whether $\{V/G\}=\Lclass^n$ in $\Kscheme{\k}$ for a
faithful $n$-dimensional $G$-representation $V$. Considering some of the
examples where we have shown that this formula is true it seems that this should
be related to the rationality of the variety $V/G$. To see if there is a more
precise relation we can look at the graded ring $\gr^*\Kscheme{\k}$ associated to
the dimension filtration of $\Kscheme{\k}$. If we let $\Bir_{\k}^n$ be the set
of birational equivalence classes of $n$-dimensional $\k$-varieties, then we
have a surjective map $\Z[\Bir_{\k}^n] \to \gr^n\Kscheme{\k}$ and it is a
natural question if it is always injective. Indeed, the injectivity for all $n$
of these maps is easily seen to be equivalent to the ``basic question'' (cf.,
\cite[Q.\ 1.2]{larsen03::motiv}) of Larsen and Lunts. If this is so, then the
relation $\{V/G\}=\Lclass^n$ in $\Kscheme{\k}$ does indeed imply that $V/G$ is
rational. Similarly, the graded ring associated to $\Kscheme{\k}'$ is equal to
the graded ring associated to $\Kscheme{\k}[\Lclass^{-1}]$ and is isomorphic to
$(\gr^0\Kscheme{\k}')[\Lclass^{-1}]$. If $\SBir_{\k}$ is equal to the stable
equivalence classes of varieties (recall that two varieties $X$ and $Y$ are
\Definition{stably birational} if $X\times\A^m$ is birational to $Y\times\A^n$
for suitable $m$ and $n$) then we have a surjection $\Z[\SBir_{\k}] \to
\gr^0\Kscheme{\k}'$ given by $[X] \mapsto \{X\}\Lclass^{-\dim X}$. If it is a
bijection then we can conclude that if $\{V/G\}=\Lclass^n$ implies that $V/G$ is
stably birational to $\A^n$.

We do not know anything about the injectivity of these maps. However, we can
make one observation. For a given $n$ we may define $\Kscheme{\k}^{\le n}$ to be
the group generated by algebraic spaces of dimension $\le n$ and with the same
relations as for $\Kscheme{\k}$ but only involving spaces of dimension $\le
n$. Then the map $\Z[\Bir_{\k}^n] \to \gr^n\Kscheme{\k}^{\le n}$ is indeed an
isomorphism. It can be checked that our calculations for $\{V/G\}$ for various
$V$ and $G$ can actually be performed in $\Kscheme{\k}^{\le n}$, where $n=\dim
V$, so that they do indeed also show that $V/G$ is rational.
\end{section}
\bibliography{preamble,abbrevs,alggeom,algebra,ekedahl}
\bibliographystyle{pretex}
\end{document}